\documentclass[12pt]{article}
\usepackage{amssymb,amsmath}
\begin{document}
\def\ddd{\displaystyle}
\def\C{{\mathbb C}}
\def\N{{\mathbb N}}
\def\Z{{\mathbb Z}}
\def\K{{\mathbb K}}
\def\R{{\mathbb R}}
\def\T{{\mathbb T}}
\def\I{{\mathbb I}}
\def\D{{\mathbb D}}
\def\zp{\N_0}
\def\M{{\cal M}}
\def\P{{\cal P}}
\def\H{{\cal H}}
\def\B{{\cal B}}
\def\pac{{\P^{\rm ac}}}
\def\pfin{{\P^{\rm fin}}}
\font\SY=msam10
\def\square{\hbox{\SY\char03}}
\def\epsilon{\varepsilon}
\def\mmu{\widehat\mu}
\def\nnu{\widehat\nu}
\def\phi{\varphi}
\def\kappa{\varkappa}
\def\theorem#1{\smallskip\noindent
{\scshape Theorem} {\bf #1}{\bf .}\hskip 8pt\sl}
\def\defin#1{\smallskip\noindent
{\scshape Definition} {\bf #1}{\bf .}\hskip 8pt\rm}
\def\prop#1{\smallskip\noindent
{\scshape Proposition} {\bf #1}{\bf .}\hskip 8pt\sl}
\def\lemma#1{\smallskip\noindent
{\scshape Lemma} {\bf #1}{\bf .}\hskip 6pt\sl}
\def\cor#1{\smallskip\noindent
{\scshape Corollary} {\bf #1}{\bf .}\hskip 6pt\sl}
\def\epr{\smallskip\rm}
\def\rem#1{\smallskip\noindent
{\scshape Remark}  {\bf #1}.\hskip 6pt}
\def\wz{\thinspace}
\def\proof{P\wz r\wz o\wz o\wz f.\hskip 6pt}
\def\quest#1{\smallskip\noindent\hskip5pt {\scshape  Question}
{\bf #1}.\hskip 6pt\sl}
\def\problem#1{\smallskip\noindent\hskip5pt {\scshape  Problem}
{\bf #1}.\hskip 6pt\sl}
\def\leq{\leqslant}
\def\geq{\geqslant}
\def\re{\text{\rm Re}\,}
\def\ker{\text{\rm ker}\,}
\def\slim{\mathop{\hbox{$\overline{\hbox{\rm lim}}$}}}
\def\ilim{\mathop{\hbox{$\underline{\hbox{\rm lim}}$}}}
\def\supp{{\rm supp}\,}
\def\ssub#1#2{#1_{{}_{{\scriptstyle #2}}}}
\def\mut{\ssub{\mu}{T}}
\def\tos{\mathrel{\hbox{$\vcenter{\offinterlineskip\halign
{\kern2pt\hfil##\hfil\kern2pt\cr
$\scriptstyle\sigma$\cr
\vrule height6pt width0pt depth0pt\cr
\smash{${-}\!\!\!\!{\longrightarrow}$}\cr
\vrule height5pt width0pt depth0pt\cr
}}$}}}
\def\ramka#1{\hbox{$\vcenter{\offinterlineskip\halign
{\vrule\vrule\kern4pt\hfil##\hfil\kern4pt\vrule\vrule\cr
\noalign{\hrule} \noalign{\hrule} \vrule height 4pt depth0pt
width0pt\cr #1\cr \vrule height 4pt depth0pt width0pt\cr
\noalign{\hrule} \noalign{\hrule}} }$}}

\title{Non-sequential weak supercyclicity and hypercyclicity}

\author{Stanislav Shkarin}

\date{}

\maketitle

\smallskip

\leftline{King's College London, Department of Mathematics}
\leftline{Strand, London WC2R 2LS, UK} \leftline{\bf e-mail: \tt
stanislav.shkarin@kcl.ac.uk}

\begin{abstract}
A bounded linear operator $T$ acting on a Banach space $\B$ is
called weakly hypercyclic if there exists $x\in \B$ such that the
orbit $\{T^n x: n=0,1,\ldots \}$ is weakly dense in $\B$ and $T$ is
called weakly supercyclic if there is $x\in \B$ for which the
projective orbit $\{\lambda T^n x: \lambda \in \C,\ n=0,1,\ldots \}$
is weakly dense in $\B$. If weak density is replaced by weak
sequential density, then $T$ is said to be weakly sequentially
hypercyclic or supercyclic respectively. It is shown that on a
separable Hilbert space there are weakly supercyclic operators which
are not weakly sequentially supercyclic. This is achieved by
constructing  a Borel probability measure $\mu$ on the unit circle
for which the Fourier coefficients vanish at infinity and the
multiplication operator $Mf(z)=zf(z)$ acting on $L_2(\mu)$ is weakly
supercyclic. It is not weakly sequentially supercyclic, since the
projective orbit under $M$ of each element in $L_2(\mu)$ is weakly
sequentially closed. This answers a question posed by Bayart and
Matheron. It is proved that the bilateral shift on $\ell_p(\Z)$,
$1\leq p <\infty$,  is weakly supercyclic  if and only if
$2<p<\infty$ and that any weakly supercyclic weighted bilateral
shift on $\ell_p(\Z)$ for $1\leq p\leq 2$ is norm supercyclic. It is
also shown that any weakly hypercyclic weighted bilateral shift on
$\ell_p(\Z)$ for $1\leq p<2$ is norm hypercyclic, which answers a
question of Chan and Sanders.
\end{abstract}

\section{Introduction}

As usual $\C$ and $\R$ are the fields of complex and real numbers
respectively, $\Z$ is the set of integers, $\N$ is the set of
positive integers and $\zp =\N\cup\{0\}$.

Let $T$ be a bounded linear operator acting on a complex Banach
space $\B$. An element $x\in\B$ is called a {\it weakly hypercyclic
vector for} $T$ if the orbit
$$
O(T,x)=\{T^nx:n\in\zp\}
$$
is weakly dense in $\B$ and $T$ is said to be {\it weakly
hypercyclic} if it has a weakly hypercyclic vector. Similarly
$x\in\B$ is called a {\it weakly supercyclic vector for} $T$ if the
projective orbit
$$
O_{\rm pr}(T,x)=\{\lambda T^nx:n\in\zp  ,\ \lambda\in\C\}
$$
is weakly dense in $\B$ and $T$ is said to be {\it weakly
supercyclic} if it has a weakly supercyclic vector.

These classes of operators are more general than the classes of
hypercyclic and supercyclic operators, in which the density is
required with respect to the norm topology, see the surveys
\cite{msa} and \cite{ge} and references therein and
\cite{cs,pra,dt,san,san1} for other related results on weak
hypercyclicity and supercyclicity. Weakly supercyclic and weakly
hypercyclic operators, although more general than the supercyclic
and hypercyclic ones, enjoy many of the properties of supercyclic
and hypercyclic operators. For instance, if $T$ is weakly
supercyclic or hypercyclic, then so is $T^n$ for any $n\in\N$. The
norm topology version of the latter result was proved by Ansari
\cite{ansa} and the same proof works for weakly supercyclic and
hypercyclic operators. Another instance: the operator $\alpha I
\oplus T: \C \oplus \B \to \C\oplus \B$, where $\B$ is a Banach
space and $\alpha \not=0$, is supercyclic if an only if
$\alpha^{-1}T$ is hypercyclic, see \cite{glm}. Again, the proof also
works if the norm topology is replaced by the weak one, see
\cite{msa} and \cite{san}.  This observation provides the first
known examples of weakly supercyclic non-supercyclic operators on a
Hilbert space \cite{san}.

Recall that a subset $A$ of a topological space $X$ is called {\it
sequentially closed} if for any convergent in $X$ sequence of
elements of $A$, the limit belongs to $A$. The minimal sequentially
closed set $[A]_{\rm s}$ containing a given set $A$ (=the
intersection of all sequentially closed sets, containing $A$) is
called the {\it sequential closure} of $A$. Finally $A\subset X$ is
called sequentially dense in $X$ if $[A]_{\rm s}=X$. Note that in
general $[A]_{\rm s}$ may be bigger than the set of limits of
converging sequences of elements of $A$.

An interesting example in the Hilbert space setting was recently
provided by Bayart and Matheron \cite{bm}. They proved that if $\mu$
is a continuous Borel probability measure on the unit circle $\T =
\{z \in \C: |z| =1\}$, supported on a Kronecker compact set, then
the multiplication operator $Mf(z)=z f(z)$ acting on $L^2(\mu)$ is
weakly supercyclic. On the other hand, since $M$ is an isometry, it
cannot be supercyclic, see \cite{ab}. It should be noted that in the
last example there is  $x$ in $\B$ such that any vector from
$L_2(\mu)$ is a limit of a weakly convergent sequence of elements of
$O_{\rm pr}(M,x)$.

The last observation motivates the following definitions. A vector
$x\in \B$ is called {\it a weakly sequentially hypercyclic} vector
for $T$ if the orbit $O(T,x)$ is weakly sequentially dense in $\B$
and $T$ is called {\it weakly sequentially hypercyclic} if it has
weakly sequentially hypercyclic vectors. A vector $x\in \B$ is
called {\it a weakly sequentially supercyclic} vector for $T$ if the
projective orbit $O_{\rm pr}(T,x)$ is weakly sequentially dense in
$\B$ and $T$ is called {\it weakly sequentially supercyclic} if it
has weakly sequentially supercyclic vectors.

Slightly different concepts were introduced by Bes, Chan and Sanders
\cite{bcs} and implicitly by Bayart and Matheron \cite{bm}. In fact,
they call the following properties weak sequential hypercyclicity
and weak sequential supercyclicity. We call them in a bit different
way in order to distinguish from the above defined ones. Namely, $T$
is called {\it weakly $1$-sequentially hypercyclic} if there exists
$x\in \B$ such that any vector from $\B$ is a limit of a weakly
convergent sequence of elements of the orbit $O(T,x)$ and $T$ is
called {\it weakly $1$-sequentially supercyclic} if there exists
$x\in \B$ such that any vector from $\B$ is a limit of a weakly
convergent sequence of elements of the projective orbit $O_{\rm
pr}(T,x)$.

The obvious relations between the above properties are summarized in
the following diagram:
$$
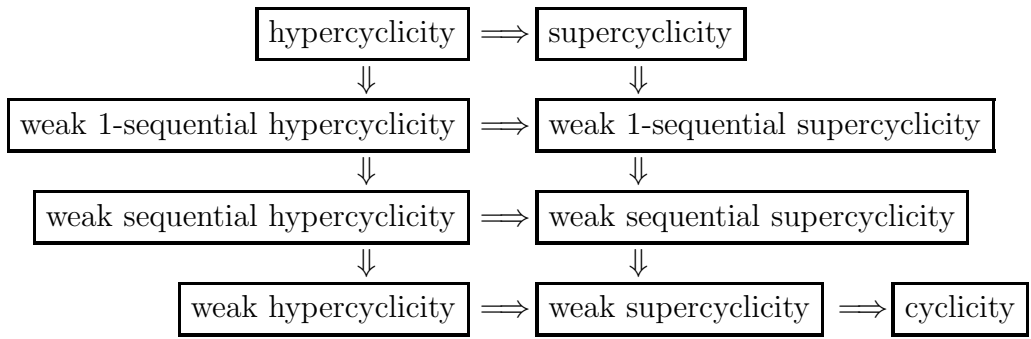
\begin{array}{rcl}
\ramka{hypercyclicity}&\!\!\!\Longrightarrow\!\!\!\!&\ramka{supercyclicity}\\
\Downarrow\qquad\quad&&\qquad\quad\Downarrow\\
\ramka{weak\ 1-sequential
hypercyclicity}&\!\!\!\Longrightarrow\!\!\!\!&\ramka{weak
1-sequential supercyclicity}\\
\Downarrow\qquad\quad&&\quad\qquad\Downarrow\\
\ramka{weak\ sequential
hypercyclicity}&\!\!\!\Longrightarrow\!\!\!\!&\ramka{weak
sequential supercyclicity}\\
\Downarrow\quad\qquad&&\quad\qquad\Downarrow\\
\ramka{weak
hypercyclicity}&\!\!\!\Longrightarrow\!\!\!\!&\ramka{weak
supercyclicity}\,\Longrightarrow\!\ramka{cyclicity}
\end{array}
$$

Bayart and Matheron \cite{bm} raised the two following questions.

\quest1 Does there exist a bounded linear operator, which is weakly
supercyclic and not weakly $1$-sequentially supercyclic? \epr

\quest2 Does there exist a positive Borel measure $\mu$ on $\T$ such
that the  Fourier coefficients of $\mu$ vanish at infinity and the
operator $Mf(z)=zf(z)$ acting on $L_2(\mu)$ is weakly supercyclic?
\epr

In view of the following proposition, an affirmative answer to the
second question implies an affirmative answer to the first one in
the Hilbert space setting.

\prop{1.1} Let $\mu$ be a non-negative Borel measure on $\T$ such
that  its Fourier coefficients $\mmu(n)=\int z^n\mu(dz)$ $(n\in\Z)$
vanish at infinity, that is $\mmu(n)\to 0$ as $|n|\to\infty$. Then
the projective orbit $O_{\rm pr}(M,f)$ is weakly sequentially closed
for any $f\in L_2(\mu)$, where the multiplication operator
$Mf(z)=zf(z)$ acts on $L_2(\mu)$. In particular, $M$ is not weakly
sequentially supercyclic. \epr

We provide an affirmative answer to Question~2 and consequently to
Question~1.

\theorem{1.2}There exists a Borel probability measure $\mu$ on $\T$
such that its  Fourier coefficients vanish at infinity and the
operator $Mf(z)=zf(z)$ acting on $L_2(\mu)$ is weakly supercyclic.
\epr

Proposition~1.1 and Theorem~1.2 immediately imply the following corollary.

\cor{1.3}There exists a weakly supercyclic unitary operator on a
separable Hilbert space, which is not weakly sequentially
supercyclic. \epr

The proof of Theorem~1.2 requires a construction  of a rather
complicated singular continuous measure. Curiously enough, it is
much easier to give an affirmative answer to Question~1 for Banach
space operators.

Given a bounded sequence $\{w_n\}_{n\in\Z}$ in $\C\setminus\{0\}$,
the {\it weighted bilateral shift} $T$ acting on $\ell_p(\Z)$,
$1\leq p<\infty$ or $c_0(\Z)$ is defined on the canonical basis
$\{e_n\}_{n\in\Z}$ by $Te_n=w_{n}e_{n-1}$. We denote
\begin{equation}
\beta (k,n)=\smash{\prod_{j=k}^n}\,|w_j| , \quad \ \ \text{for}\ \
k,n\in\Z\ \text { with } \ k\leq n. \label{beta}
\end{equation}
In the particular case $w_n\equiv1$ we have the {\it unweighted
bilateral shift}, which we denote as $B$.

Salas \cite{salas1,salas} has characterized hypercyclic and
supercyclic bilateral weighted shifts in terms of weight sequences.
We formulate his results in a slightly different form, however
obviously equivalent to the original ones.

\theorem{S}Let $T$ be a bilateral weighted shift acting on
$\ell_p(\Z)$ with $1\leq p<\infty$ or $c_0(\Z)$. Then $T$ is
hypercyclic if and only if for any $k\in\zp$,
\begin{equation}
\ilim\limits_{n\to\infty} \max\biggl\{\max\limits_{|j|\leq
k}\beta(j-n,j), \Bigl(\min\limits_{|j|\leq
k}\beta(j,j+n)\Bigr)^{-1}\biggr\}=0 \label{sal1}
\end{equation} \epr
and $T$ is supercyclic if and only if for any $k\in\zp$,
\begin{equation}
\ilim\limits_{n\to+\infty} \Bigl(\max\limits_{|j|\leq k}
\beta(j-n,j)\Bigr)\Bigl(\min\limits_{|j| \leq
k}\beta(j,j+n)\Bigr)^{-1}=0. \label{sal2}
\end{equation}\rm

This theorem implies, in particular, that hypercyclicity and
supercyclicity of a bilateral weighted shift acting on $\ell_p(\Z)$
with $1\leq p<\infty$ do not depend on $p$. It will be clear from
the results below that it is not the case for weak hypercyclicity
and weak supercyclicity.

The main result of the paper \cite{bcs} by Bes, Chan and Sanders is
the following.

\theorem{BCS}Let $T$ be a bilateral weighted shift acting on
$\ell_p(\Z)$, $1\leq p<\infty$. If $T$ is weakly $1$-sequentially
hypercyclic then $T$ is hypercyclic. If $T$ is weakly
$1$-sequentially supercyclic then $T$ is supercyclic. \epr

We prove the following slightly stronger statement.

\prop{1.4}Let $T$ be a bilateral weighted shift acting on
$\ell_p(\Z)$, $1\leq p<\infty$ or $c_0(\Z)$. If $T$ is weakly
sequentially hypercyclic then $T$ is hypercyclic. If $T$ is weakly
sequentially supercyclic then $T$ is supercyclic. \epr

In \cite{san1} it is proved that the unweighted bilateral shift $B$
acting on $c_0(\Z)$ is weakly supercyclic. This result is a
corollary of the following stronger one.

\theorem{1.5}The unweighted bilateral shift $B$ on $\ell_p(\Z)$ is
weakly supercyclic if and only if $p>2$. \epr

Thus, $B$ acting on $\ell_p(\Z)$ for $2<p<\infty$ provides an
example of a weakly supercyclic not weakly sequentially supercyclic
isometric linear operator acting on a uniformly convex Banach space.
Since any $\ell_p(\Z)$ is densely and continuously embedded into
$c_0(\Z)$, Theorem~1.5, via comparison principle, implies weak
supercyclicity of $B$ on $c_0(\Z)$. It worth mentioning that the
proof of the above result is completely different from the one in
\cite{san1} for $c_0(\Z)$. Theorems~1.2 and~1.5 are in strong
contrast with Ansari and Bourdon's result \cite{ab} that a Banach
space isometry can not be supercyclic. In \cite{cs} Chan and Sanders
have shown that

\theorem{CS}The bilateral weighted shift $T$ with the weight
sequence $w_n=2$ if $n\geq 0$, $w_n=1$ if $n<0$ acting on
$\ell_p(\Z)$ for $2\leq p<\infty$ is weakly hypercyclic and
non-hypercyclic. \epr

They also raised the following natural question.

\quest3 Does there exist a weakly hypercyclic non-hypercyclic
bilateral weighted shift acting on $\ell_p(\Z)$ for $1\leq p<2$?
\epr

We answer this question negatively.

\theorem{1.6} Let $T$ be a bilateral weighted shift acting on
$\ell_p(\Z)$. If $1\leq p<2$ and $T$ is weakly hypercyclic then $T$
is hypercyclic. If $1\leq p\leq 2$ and $T$ is weakly supercyclic
then $T$ is supercyclic. \rm

Theorem~1.6,  Proposition~1.4 and Theorem~CS immediately imply the
following corollary.

\cor{1.7} Let $1\leq p<\infty$. Then any weakly hypercyclic
bilateral weighted shift acting on $\ell_p(\Z)$ is hypercyclic if
and only if $p<2$. Moreover any weakly supercyclic bilateral
weighted shift acting on $\ell_p(\Z)$ is supercyclic if and only if
$p\leq2$. \epr

Bes, Chan and Sanders \cite{bcs} have also raised the following
questions.

\quest4 Does there exist an invertible bilateral weighted shift $T$
acting on $\ell_p(\Z)$ such that $T$ and $T^{-1}$ are both weakly
hypercyclic and $T$ is not hypercyclic? Does there exist a weakly
hypercyclic bilateral weighted shift $T$, acting on $\ell_p(\Z)$
such that $T$ is not supercyclic? \epr

We answer both questions affirmatively:

\prop{1.8}There exists an invertible non-hypercyclic bilateral
weighted shift $T$ acting on $\ell_2(\Z)$ such that both $T$ and
$T^{-1}$ are weakly hypercyclic. \epr

\prop{1.9}For any $p>2$ there exists a weakly hypercyclic bilateral
weighted shift acting on $\ell_p(\Z)$, which is not supercyclic.
\epr

Section~2 is devoted to some basic facts about weak limit points,
their relation with $p$-sequences and antisupercyclicity, a concept
that was introduced in \cite{shk}. Proposition~1.4 is proved in the
end of Section~2. In Section~3 Theorem~1.4, Proposition~1.8 and
Proposition~1.9 are proved. In Section~4 we prove Proposition~1.1
and Theorem~1.2, which is probably the most difficult result in this
article. Theorem~1.6 is proved in Section~5. In Section~6 we discuss
the tightness of certain results of the previous sections and pose
some open questions related to this work.

\section{Antisupercyclicity and weak closures}

Throughout this section $\Lambda$ is an infinite countable set.
Recall that $\ell_\infty(\Lambda)$ is the space of complex valued or
real valued bounded sequences $\{x_\alpha\}_{\alpha \in \Lambda }$
endowed with the supremum norm and $c_0(\Lambda)$ is the subspace of
$\ell_\infty(\Lambda)$ consisting of  sequences
$\{x_\alpha\}_{\alpha\in\Lambda}$  such that $\{\alpha\in\Lambda:
|x_\alpha|>\epsilon\}$ is finite for each $\epsilon
>0$. For $1\leq p <\infty$, $\ell_p(\Lambda)$ is the space
of sequences  $x\in\ell_\infty(\Lambda)$ for which
$$
\|x\|_p=\biggl(\sum_{\alpha\in\Lambda }|x_{\alpha}|^p\biggr)^{1/p}
<\infty.
$$
Of course, these spaces are isomorphic to the usual sequence spaces
$\ell_p$ and $c_0$ indexed on $\zp$. The point is that sometimes it
is more convenient to specify a different index set. For each
$\alpha\in\Lambda$, we denote by $e_\alpha$ the sequence in which
all elements, except the $\alpha$-th, whose value is one, vanish. It
is well-known that $\{e_\alpha\}_{\alpha\in\Lambda }$ is an
unconditional absolute Schauder basis in $\ell_p(\Lambda)$ for
$1\leq p <\infty$ and in $c_0(\Lambda)$. This basis is usually
called {\it the canonical basis}. For $x\in\ell_p(\Lambda)$ and
$y\in\ell_q(\Lambda)$ with $\frac1p+\frac1q=1$, we denote
$$
\langle x,y\rangle=\sum_{\alpha\in\Lambda}x_\alpha y_\alpha.
$$
In what follows for a sequence $x=\{x_\alpha\}_{\alpha\in\Lambda}$
we shall usually write $\langle x,e_\alpha\rangle$ instead of
$x_\alpha$. The {\it support} of a sequence
$x=\{x_\alpha\}_{\alpha\in\Lambda}$ is the set
$$
\supp(x)=\{\alpha\in\Lambda:x_\alpha\neq0\}=\{\alpha\in\Lambda:
\langle x,e_\alpha\rangle\neq0\}.
$$

\subsection{Antisupercyclicity}

A bounded linear operator $T$ acting on a Banach space $\B$ is
called {\it antisupercyclic} if the sequence $\{T^n
x/\|T^nx\|\}_{n\in\zp}$ converges weakly to zero for any $x\in \B$
such that  $T^n x\neq 0$ for each $n\in\zp$.  This is the case when
the angle criterion of supercyclicity \cite{msa} is not satisfied in
the strongest possible way. In Hilbert space antisupercyclicity
means that the angles between any fixed vector $y$ and the elements
$T^nx$ of any orbit, not vanishing eventually, tend to $\pi/2$.

\theorem{2.1}Let $T$ be an antisupercyclic bounded  linear operator
acting on a Banach space $\B$. Then for any $x\in \B$, the
projective orbit $O_{\rm pr}(T,x)$ is weakly sequentially closed in
$\B$. In particular, antisupercyclic operators are never weakly
sequentially supercyclic if $\dim \B>1$. \epr

\proof Let $x\in\B$ and $\{y_n\}$ be a weakly convergent sequence of
elements of $O_{\rm pr}(T,x)$. For any $m\in\zp$, let $L_m=\{\lambda
T^nx: \lambda\in\C,\ 0\leq n\leq m\}$. If each $y_n$ belongs to
$L_m$ for some $m$, then taking into account that $L_m$ is weakly
closed, we see that the weak limit of the sequence $y_n$ belongs to
$L_m\subset O_{\rm pr}(T,x)$. Otherwise, $T^nx\neq 0$ for any
$n\in\zp$ and passing to a subsequence, if necessary, we can assume
that $y_n=(c_n/\|T^{m_n}x\|)T^{m_n}x$, where $c_n\in\C$ and $m_n$ is
a strictly increasing sequence of positive integers. Since any
weakly convergent sequence is bounded, we find that $\{c_n\}$ is
bounded. Antisupercyclicity of $T$ implies that
$z_n=T^{m_n}x/\|T^{m_n}x\|$ tends weakly to zero. Since $c_n$ is
bounded, we conclude that $y_n=c_nz_n$ tends weakly to zero, which
is in $O_{\rm pr}(T,x)$. Hence $O_{\rm pr}(T,x)$ is weakly
sequentially closed. \square

\subsection{Weak density and $p$-sequences}

\lemma{2.2}Let $1<p\leq\infty$, $\{c_\alpha\}_{\alpha\in\Lambda}$ be
a sequence of complex numbers and $\B_p=\ell_p(\Lambda)$ if
$1<p<\infty$, $\B_\infty=c_0(\Lambda)$. Then zero is in the weak
closure of the set $Y=\{c_\alpha e_\alpha:\alpha\in\Lambda\}$ in the
Banach space $\B_p$ if and only if
\begin{equation}
\sum_{\alpha\in\Lambda}|c_\alpha|^{-q}=\infty,\ \ \text{where \
$\textstyle \frac1q+\frac1p=1$}. \label{l24}
\end{equation}
\rm

\proof Without loss of generality, we may assume that
$c_\alpha\neq0$ for each $\alpha\in\Lambda$, otherwise the result is
trivial. Assume that (\ref{l24}) is not satisfied. Then $b\in
\ell_q(\Lambda)=\B_p^*$, where $b_\alpha=|c_\alpha|^{-1}$,
$\alpha\in\Lambda$. Clearly $|\langle c_\alpha e_\alpha,b\rangle|=1$
for any $\alpha\in\Lambda$. Therefore zero is not in the weak
closure of $Y$.

Conversely assume that (\ref{l24}) is satisfied. Then $b\notin
\ell_q(\Lambda)$. Let $x_1,\dots,x_m\in\ell_q(\Lambda)=\B_p^*$ and
$a_\alpha=\max\limits_{1\leq j\leq m} |\langle
x_j,e_\alpha\rangle|$. Since $a=\{a_\alpha\}\in\ell_q(\Lambda)$,
$b\notin\ell_q(\Lambda)$ and the entries of $a$ and $b$ are
non-negative, we have $\inf\limits_{\alpha\in\Lambda} a_\alpha
b_\alpha^{-1}=0$. Finally observe that  $|\langle c_\alpha
e_\alpha,x_j\rangle|\leq a_\alpha b_\alpha^{-1}$ for any
$\alpha\in\Lambda$ and $1\leq j\leq m$. Hence
$\inf\limits_{\alpha\in\Lambda}\max\limits_{1\leq j\leq m} |\langle
c_\alpha e_\alpha,x_j\rangle|=0$. Thus, zero is in the weak closure
of $Y$. \square

Let $1\leq p\leq\infty$. A sequence
$\{x_\alpha\}_{\alpha\in\Lambda}$ of elements of a Banach space $\B$
is called a $p$-{\it sequence} if there exists $c>0$ such that
\begin{equation}
\left\|\sum_{j=1}^n a_jx_{\alpha_j}\right\|\leq c\|a\|_p\ \
\text{for any $n\in\N$, any $a\in\C^n$ and any pairwise different
$\alpha_1,\dots,\alpha_n\in\Lambda$}. \label{pseq}
\end{equation}

For instance, each bounded sequence in $\ell_p$ with disjoint
supports is a $p$-sequence. Clearly $\{x_\alpha\}$ is a $p$-{\it
sequence} if and only if there exists a bounded linear operator
$S:\B_p\to \B$ such that $Se_\alpha=x_\alpha$ for each
$\alpha\in\Lambda$, where $\B_p=\ell_p(\Lambda)$ if $1\leq
p<\infty$, $\B_\infty=c_0(\Lambda)$. The concept of $p$-sequence
provides an easy sufficient condition for zero to belong to the weak
closure of certain sequences.

\lemma{2.3}Let $1<p\leq\infty$ and $\{x_\alpha\}_{\alpha\in\Lambda}$
be a $p$-sequence in a Banach space $\B$ and
$\{c_\alpha\}_{\alpha\in\Lambda}$ be a sequence of complex numbers,
satisfying $(\ref{l24})$. Then zero is in the weak closure of
$Y=\{c_\alpha x_\alpha:\alpha\in\Lambda\}$ in $\B$. \epr

\proof Let $\B_p=\ell_p(\Lambda)$ if $1<p<\infty$ and
$\B_\infty=c_0(\Lambda)$. Since $\{x_\alpha\}$ is a $p$-sequence,
there  exists a bounded linear operator $S:\B_p\to \B$ such that
$Se_\alpha=x_\alpha$ for each $\alpha\in\Lambda$. By Lemma~2.2, zero
is in the  weak closure of the set $N=\{c_\alpha
e_\alpha:\alpha\in\Lambda\}$  in $\B_p$. Since $S(N)=Y$ and
$S:\B_p\to \B$ is weak-to-weak continuous, we see that zero is in
the weak closure of $Y$. \square

The previous lemma allows us to prove the following proposition,
which provides sufficient conditions for weak supercyclicity and
hypercyclicity.

\prop{2.4} Let $\B$ be a Banach space, $T:\B\to\B$ be a bounded
linear operator, $S$ be a subset of $\B$ such that $\Omega=\{\lambda
x:\lambda\in\C,\ x\in S\}$ is weakly dense in $\B$ and $u\in\B$.
Assume also that for any $x\in S$, there exist $p_x\in(1,\infty]$,
an infinite set $A_x\subset\zp$ and  maps
$\alpha_x,\beta_x:A_x\to\C$ and $\gamma_x:A_x\to\N$ satisfying
\begin{itemize}\itemsep=-2pt
\item[{\rm(C1)}]$\{\beta_x(k)T^{\gamma_x(k)}u-\alpha_x(k)x\}_{k\in A_x}$ is a
$p_x$-sequence in $\B$;
\item[{\rm(C2)}]$\sum\limits_{k\in A_x}|\alpha_x(k)|^{q_x}=\infty$, where
$\frac1{p_x}+\frac1{q_x}=1$.
\end{itemize}
Then $u$ is a weakly supercyclic vector for $T$.

If additionally $S$ itself is weakly dense in $\B$ and
$\alpha_x=\beta_x$ for each $x\in S$, then $u$ is a weakly
hypercyclic vector for $T$. \epr

\proof Let $x\in S$. Lemma~2.3 along with (C1) and (C2) implies that
zero is in the weak closure of
$\{\frac{\beta_x(k)}{\alpha_x(k)}T^{\gamma_x(k)}u-x:k\in A_x\}$.
Thus, $x$ is in the weak closure of
$\{\frac{\beta_x(k)}{\alpha_x(k)}T^{\gamma_x(k)}u:k\in A_x\}$, which
is contained in $O_{\rm pr}(T,u)$. Since $x$ is an arbitrary element
of $S$ and $O_{\rm pr}(T,u)$ is stable under the multiplication by
scalars, we see that $\Omega$ is contained in the weak closure of
$O_{\rm pr}(T,x)$. Taking into account that $\Omega$ is weakly dense
in $\B$, we observe that $O_{\rm pr}(T,u)$ is weakly dense in $\B$.
Thus, $u$ is a weakly supercyclic vector for $T$. Suppose now that
$S$ is weakly dense in $\B$ and $\alpha_x=\beta_x$ for each $x\in
S$. Then any $x\in S$ is in the weak closure of
$\{T^{\gamma_x(k)}u:k\in A_x\}\subseteq O(T,u)$. Therefore $O(T,u)$
is weakly dense in $\B$. Thus, $u$ is a weakly hypercyclic vector
for $T$. \square

The following lemma deals with perturbations of $p$-sequences.

\lemma{2.5}Let $\{x_\alpha\}_{\alpha\in\Lambda}$ and
$\{y_\alpha\}_{\alpha\in\Lambda}$ be two sequences in a Banach space
$\B$, where the first one is a $p$-sequence for $1\leq p\leq\infty$
and $b\in \ell_q(\Lambda)$, where  $b_\alpha=\|x_\alpha-y_\alpha\|$
and $\frac1p+\frac1q=1$. Then $\{y_\alpha\}_{\alpha\in\Lambda}$ is a
$p$-sequence. \epr

\proof Since $\{x_\alpha\}_{\alpha\in\Lambda}$  is a $p$-sequence,
there exists $c>0$ such that (\ref{pseq}) is satisfied. Let
$n\in\N$, $a\in\C^n$ and $\alpha_1,\dots,\alpha_n$ be pairwise
different elements of $\Lambda$. Upon applying the H\"older
inequality, we obtain
\begin{equation*}
\left\|\sum_{j=1}^na_jy_{\alpha_j}\right\|\leq
\left\|\sum_{j=1}^na_jx_{\alpha_j}\right\|+
\left\|\sum_{j=1}^na_j(x_{\alpha_j}-y_{\alpha_j}) \right\|\leq
c\|a\|_p+\sum_{j=1}^n|a_j|b_{\alpha_j} \leq (c+\|b\|_q) \|a\|_p.
\end{equation*}
Hence $\{y_\alpha\}_{\alpha\in\Lambda}$ is a $p$-sequence. \square

We end this section with a sufficient condition for being a
2-sequence in a Hilbert space.

\lemma{2.6}Let $\{g_{n}\}_{n\in\N}$ be a bounded sequence in a
Hilbert space $\H$ such that
$$
c=\sum_{1\leq m<n<\infty}|\langle g_n,g_m\rangle|^2<\infty.
$$
Then $\{g_n\}_{n\in\N}$ is a $2$-sequence in $\H$. \epr

\proof Denote $d=\sup\limits_{n\in\N}\|g_n\|$ and let $n\in\N$,
$a\in\C^n$ and $m_1,\dots,m_n$ be pairwise different positive
integers. Applying the Cauchy--Schwartz inequality, we obtain
\begin{align*}
\left\|\sum_{j=1}^na_jg_{m_j}\right\|^2&=\left\langle
\sum_{j=1}^na_jg_{m_j},\sum_{k=1}^na_kg_{m_k}\right\rangle=
\sum_{j=1}^n\sum_{k=1}^n a_j\overline{a_k} \langle
g_{m_j},g_{m_k}\rangle\leq\sum_{j=1}^n|a_j|^2\|g_{m_j}\|^2+
\\
&\quad+\!\!\! \sum_{1\leq j,k\leq n\atop j\neq k} |a_ja_k\langle
g_{m_j},g_{m_k}\rangle|\leq d^2\|a\|_2^2+ \biggl(\sum_{1\leq j,k\leq
n\atop j\neq k} |a_ja_k|^2\biggr)^{1/2}\!\biggl(\sum_{1\leq j,k\leq
n\atop j\neq k} |\langle
g_{m_j},g_{m_k}\rangle|^2\biggr)^{1/2}\!\!\!\leq
\\
&\leq d^2\|a\|_2^2+ \biggl(\sum_{1\leq j,k\leq n}
|a_ja_k|^2\biggr)^{1/2}\!\biggl(2\sum_{j,k\in\N\atop j>k} |\langle
g_j,g_k\rangle|^2\biggr)^{1/2}=(d^2+(2c)^{1/2})\|a\|_2^2.
\end{align*}
Hence (\ref{pseq}) for $\Lambda=\N$, $x_n=g_n$ and $p=2$ is
satisfied with the constant $(d^2+\sqrt{2c})^{1/2}$. Thus,
$\{g_n\}_{n\in\N}$ is a 2-sequence. \square

\subsection{Proof of Proposition 1.4}

In \cite{shk} it is proven that

\theorem{A}A weighted bilateral shift $T$ acting on $\ell_p(\Z)$,
$1<p<\infty$ is antisupercyclic if and only if it is not
supercyclic. \epr

The same result is true and the same proof works when $T$ is acting
on $c_0(\Z)$. One has to take into account that $c_0(\Lambda)$
shares the following property with $\ell_p(\Lambda)$ for
$1<p<\infty$: a sequence is weakly convergent if and only if it is
norm bounded and coordinatewise convergent. This fails for sequences
in $\ell_1(\Lambda)$ and so does the above theorem.

Let $\B_p=\ell_p(\Z)$ if $1<p<\infty$ and $\B_\infty=c_0(\Z)$. The
above result along with Proposition~2.1 and the comparison principle
implies that a weighted bilateral shift $T$ acting on $\B_p$ for
$1\leq p\leq\infty$ is weakly sequentially supercyclic if and only
if it is supercyclic. Moreover the projective orbits of $T$ are
weakly sequentially closed if $T$ is not supercyclic. It remains to
show that a weakly sequentially hypercyclic bilateral weighted shift
is hypercyclic. The situation with hypercyclicity differs from that
with supercyclicity since, as it was mentioned in \cite{bcs}, orbits
of a non-hypercyclic weighted bilateral shift may be not weakly
sequentially closed. The proof goes along the same lines as in
\cite{bcs}, but we have to overcome few additional difficulties.

Let $1\leq p\leq\infty$ and $T$ be a weakly sequentially hypercyclic
bilateral weighted shift on $\B_p$. Denote by $\Omega_0$ the set of
weakly hypercyclic vectors for $T$. Let $\Omega$ be the set of
$z\in\B_p$ for which there exist $x\in\Omega$ and a strictly
increasing sequence $\{n_k\}_{k\in\zp}$ of positive integers such
that the sequence $T^{n_k}x$ is weakly convergent to $z$. Since $T$
is weakly sequentially hypercyclic, $\Omega$ is weakly sequentially
dense in $\B_p$.

\lemma{2.7}For any $z\in\Omega$, any $l\in\N$, any $\epsilon>0$ and
any $y\in\B_p$ with finite support, there exists $v\in\B_p$ with
finite support and $n\in\N$ such that $n>l$, $\|v\|_p<\epsilon$,
$\|T^ny\|_p<\epsilon$ and $\|T^nv-z\|_p<\epsilon$. \epr

\proof Since $z\in\Omega$, there exist a weakly hypercyclic vector
$x$ for $T$ and a strictly increasing sequence $n_k$ of positive
integers such that $T^{n_k}x$ converges weakly to $z$ as
$k\to\infty$. Since any weakly convergent sequence is bounded, there
exists $M>0$ such that $\|T^{n_k}x\|_p\leq M$ for any $k\in\zp$.
Clearly $\|u-P_{0,d} u\|_p\to 0$ as $d\to\infty$ for any $u\in\B_p$,
where
$$
P_{a,d}:\B_p\to \B_p,\qquad P_{a,d} u=\sum_{n=a-d}^{a+d} \langle
u,e_n\rangle e_n.
$$
Pick $r\in\N$ such that $\|z-P_{0,r}z\|_p<\epsilon/2$. Since $x$ is
a weakly hypercyclic vector for $T$ and $y$ has finite support,
there exists $m\in\zp$ such that $|\langle
T^mx,e_j\rangle|>M|\langle y,e_j\rangle|/\epsilon$, whenever
$j\in\supp(y)$. Taking into account that $T$ is a weighted shift, we
see that for any $l\in\zp$,
$$
|\langle T^{m+l}x,e_j\rangle|>M|\langle T^ly,e_j\rangle|/\epsilon,\
\ \text{whenever}\ \ \langle T^ly,e_j\rangle\neq 0.
$$
Hence $\|T^ly\|_p< \frac{\epsilon}{M}\|T^{m+l}x\|_p$ for each
$l\in\zp$. In particular, for $l=n_k-m$, we have
$$
\|T^{n_k-m}y\|_p<\frac{\epsilon}{M}\|T^{n_k}x\|_p\leq \epsilon, \ \
\text{whenever}\ \ n_k\geq m.
$$
Denote $v_k=T^mP_{n_k,r}x$. Clearly $\|v_k\|_p\to 0$ as
$k\to\infty$. Since $T^{n_k-m}v_k=P_{0,r}T^{n_k}x$ and $T^{n_k}x$
converges weakly to $z$, we see that $T^{n_k-m}v_k$ converges weakly
to $P_{0,r} z$. Since all the vectors $T^{n_k-m}v_k$ belong to the
finite dimensional range of $P_{0,r}$, we have
$\|T^{n_k-m}v_k-P_{0,r} z\|_p\to 0$. Choosing $k$ large enough, we
can ensure that $n_k-m>l$, $\|v_k\|_p<\epsilon$ and
$\|T^{n_k-m}v_k-P_{0,r} z\|_p<\epsilon/2$. Since
$\|z-P_{0,r}z\|_p<\epsilon/2$, we have
$\|T^{n_k-m}v_k-z\|_p<\epsilon$. Thus, $v=v_k$ and $n=n_k-m$ satisfy
all desired conditions.  \square

\lemma{2.8}For any sequence $\{z_k\}_{k\in\zp}$ of elements of \
$\Omega$ \ there exist a weakly hypercyclic for $T$ vector
$x\in\B_p$ and a strictly increasing sequence $\{n_k\}_{k\in\zp}$ of
positive integers such that $\|T^{n_k}x-z_k\|_p\to 0$ as
$k\to\infty$. \epr

\proof Since any sequence of elements of $\Omega$ is a subsequence
of a sequence of elements of $\Omega$, which is norm-dense in
$\Omega$, we can, without loss of generality assume that
$\{z_n:n\in\zp\}$ is norm-dense in $\Omega$.

By Lemma~2.7 we can construct inductively a strictly increasing
sequence $\{n_k\}_{k\in\zp}$ of positive integers and
$\{x_k\}_{k\in\zp}$ of vectors in $\B_p$ with finite supports such
that
\begin{align*}
&\|x_k\|_p<s_k, \ \ \|T^{n_k}u_k\|_p<s_k\ \ \text{and}\ \
\|T^{n_k}x_k-z_k\|_p<s_k, \ \ \text{where}
\\
&u_0=0,\ \ u_k=x_0+{\dots}+x_{k-1}\ \ \text{if}\ \ k\geq1,\ \ s_0=1\
\ \text{and}
\\
&s_k=2^{-k}\min\{1,\|T^{n_0}\|_p^{-1},\dots,\|T^{n_{k-1}}\|_p^{-1}\}\
\ \text{if}\ \ k\geq1.
\end{align*}

Since $\|x_k\|_p\leq 2^{-k}$ the series $\sum\limits_{k=0}^\infty
x_k$ is absolutely convergent in $\ell_p(\Z)$. Let
$x=\sum\limits_{k=0}^\infty x_k$. Then
\begin{align*}
\|T^{n_k}x-z_k\|_p&=\biggl\|(T^{n_k}x_k-z_k)+T^{n_k}u_k+
\sum_{j=k+1}^\infty T^{n_k}x_j\biggr\|_p\leq
\\
&\leq 2s_k+\sum_{j=k+1}^\infty \|T^{n_k}\|_ps_j\leq
2^{-k+1}+\sum_{j=k+1}^\infty 2^{-j}=3\cdot2^{-k}\to 0\ \ \text{as
$k\to\infty$.}
\end{align*}
Since $\{z_n:n\in\zp\}$ is norm-dense in $\Omega$, $\Omega$ is
weakly dense in $\B_p$ and $\|T^{n_k}x-z_k\|_p\to0$, we see that $x$
is a weakly hypercyclic vector for $T$. Thus, $x$ and $\{n_k\}$
satisfy all desired conditions. \square

Let $\{u_n\}_{n\in\zp}$ be a sequence of  elements of $\Omega$
weakly convergent to $u\in\B_p$. By Lemma~2.8 there exist a weakly
hypercyclic for $T$ vector $x\in\B_p$ and a strictly increasing
sequence $\{k_n\}_{n\in\zp}$ of positive integers such that
$\|T^{k_n}x-u_n\|_p\to 0$ as $n\to\infty$. Since $u_n$ tends weakly
to $u$, we have that $T^{k_{n}}x$ tends weakly to $u$. Hence
$u\in\Omega$ and therefore $\Omega$ is weakly sequentially closed.
Since $\Omega$ is weakly sequentially dense in $\B_p$, we have
$\Omega=\B_p$. Taking a norm dense sequence $\{f_n\}_{n\in\zp}$ in
$\Omega=\B_p$ and applying Lemma~2.8 once again, we obtain
$y\in\B_p$ and a strictly increasing sequence $m_n$ of positive
integers such that $\|T^{m_{n}}y-f_n\|_p\to 0$. It follows that
$O(T,y)$ is norm dense in $\B_p$. Hence $y$ is a hypercyclic vector
for $T$. The proof is complete.

\section{Weakly supercyclic and hypercyclic bilateral shifts}

Before proving Theorem~1.5, and Propositions~1.9 and~1.10, we, using
Proposition~2.4, shall derive sufficient conditions for weak
hypercyclicity and weak supercyclicity of invertible weighted shifts
in terms of weight sequences. Our sufficient condition of weak
hypercyclicity of a bilateral weighted shift differs from the one of
Chan and Sanders \cite{cs} and is fairly easier to handle. In fact
it is possible, using basically the same proof, to generalize our
criteria for non-invertible bilateral weighted shifts, but the
conditions become too heavy in this case.

Recall that the {\it density} of a subset $A\subset\zp$ is the limit
$\lim\limits_{n\to\infty}\frac{N(n)}{n}$, where $N$ is the counting
function of $A$, that is, $N(n)$ is the number of elements of the
set $\{m\in A:m\leq n\}$. The following elementary lemma can be
found in many places, see for instance \cite{levin}, Chapter~1.

\lemma{3.1}Let $A$ be a subset of $\zp$ of positive density and
$\{s_n\}_{n\in\zp}$ be a monotonic sequence of positive numbers.
Then $\sum\limits_{n\in\zp}s_n=\infty$ if and only if
$\sum\limits_{n\in A}s_n=\infty$. \epr

For a sequence $x=\{x_n\}_{n\in\Z}$ of complex numbers denote
$$
\gamma(x)=\max_{n\in\supp(x)}|n|.
$$

In the following two lemmas $\B_p=\ell_p(\Z)$ if $1\leq p<\infty$
and $\B_\infty=c_0(\Z)$.

\lemma{3.2}Let $1\leq p<\infty$, and $\{a_n\}_{n\in\zp}$ be a
sequence of non-negative numbers such that
$\lim\limits_{n\to\infty}a_n=\infty$. Then there exists a map
$\kappa:\zp\to\B_p$ such that
\begin{itemize}\itemsep=-2pt
\item[{\rm (U1)}]the set $S=\kappa(\zp)$ consists of vectors with finite
support and is norm-dense in $\B_p$;
\item[{\rm (U2)}]$\gamma(\kappa(n))\leq a_n$ and
$\|\kappa(n)\|_p\leq a_n$ for each $n\in\zp$;
\item[{\rm (U3)}]for each $x\in S\setminus\{0\}$ the set $\kappa^{-1}(x)$
has positive density.
\end{itemize}\epr

\proof Take a dense in $\B_p$ sequence $\{x_n\}_{n\in\zp}$ of
pairwise different non-zero vectors with finite support. Since
$\lim\limits_{n\to\infty} a_n=\infty$, we can pick a strictly
increasing sequence $\{m_n\}_{n\in\zp}$ of positive integers such
that $\|x_n\|_p\leq a_k$ and $\gamma(x_n)\leq a_k$ for each $k\geq
m_n$. Choose a strictly increasing sequence $\{p_n\}_{n\in\zp}$ of
prime numbers such that $p_n>m_n$ for any $n\in\zp$ and put
\begin{equation*}
A_0=\{jp_0+1:j\in\N\}\ \ \text{and}\ \ A_n=\{p_0\cdot{\dots}\cdot
p_{n-1}(jp_n+1):j\in\N\}\ \ \text{for}\ \ n\geq1.
\end{equation*}
Each $A_n$ is an arithmetic progression and therefore has positive
density. From easy divisibility considerations it follows that the
sets $A_n$ are disjoint. This allows us to define the map kappa by
setting $\kappa(m)=x_n$ if $m\in A_n$ and $\kappa(m)=0$ if
$m\in\zp\setminus\bigcup\limits_{n=0}^\infty A_n$.

Since $S=\kappa(\zp)=\{0\}\cup\{x_n:n\in\zp\}$, we see that $S$ is
dense in $\ell_p(\Z)$ and consists of vectors with finite support.
Since $\kappa^{-1}(x_n)=A_n$, condition (U3) is satisfied. Let $m\in
A_n$. From the definition of $A_n$ it follows that $m\geq m_n$ and
therefore $\kappa(m)=x_n$ satisfies (U2). If
$m\in\zp\setminus\bigcup\limits_{n=0}^\infty A_n$, then
$\kappa(m)=0$ and (U2) is trivially satisfied. Thus, $\kappa$
satisfies all required conditions. \square

For $n\in\zp$ and $m\in\Z$ denote
\begin{equation}
L(m,n)=\{k\in\Z:|k-m|\leq n\}.\label{L}
\end{equation}
Clearly
\begin{equation}
\text{$L(a,b)\cap L(c,d)=\varnothing$ if and only if $|a-c|>b+d$}.
\label{ldis}
\end{equation}

\lemma{3.3}Let $T$ be a bilateral weighted shift acting on $\B_p$,
$1\leq p\leq\infty$, $\{a_n\}_{n\in\zp}$, $\{r_n\}_{n\in\zp}$ be
monotonically non-decreasing sequences of non-negative integers such
that such that $r_n-r_{n-1}-r_{n-2}> a_{n}+a_{n-1}$ for any
$n\geq2$, $\{x_{n,k}\}_{n,k\in\zp}$ be a double sequence of vectors
from $\ell_p(\Z)$ such that $\gamma(x_{n,k})\leq a_n$ for each
$n,k\in\zp$ and
$$
y_k=\sum_{n\in\zp,\ n\neq k\atop a_n<
(r_k-r_{k-1})/2}T^{r_k-r_n}x_{n,k}\in\B_p,\ \ k\in\N.
$$
Then $y_k$ have disjoint supports. \epr

\proof Since $\gamma(x_{n,k})\leq a_n$, we have
$\supp(x_{n,k})\subseteq L(0,a_n)$ for each $n,k\in\zp$, where the
sets $L(m,n)$ are defined in (\ref{L}). Therefore
$\supp(T^{r_k-r_n}x_{n,k})\subseteq L(r_n-r_k,a_n)$ for each
$n,k\in\zp$. Hence,
$$
\supp(y_k)\subseteq\bigcup_{n\in\zp,\ n\neq k\atop a_n<
(r_k-r_{k-1})/2} L(r_n-r_k,a_n).
$$
Let $k,l\in\N$ and $k>l$. We have to show that
$\supp(y_k)\cap\supp(y_l)=\varnothing$. According to the last
display and (\ref{ldis}) it suffices to verify that
\begin{equation}
\begin{array}{l}
|(r_n-r_k)-(r_m-r_l)|>a_n+a_m
\\
\text{if $m,n\in\zp$, $a_n<(r_k-r_{k-1})/2$, $a_m<(r_l-r_{l-1})/2$,
$n\neq k$ and $m\neq l$.}
\end{array}
\label{44}
\end{equation}
Let $m,n\in\zp$ be such that $a_n<(r_k-r_{k-1})/2$, $a_m<
(r_l-r_{l-1})/2$, $n\neq k$ and $m\neq l$.

{\bf Case} $m\neq n$. Denote $j=\max\{n,k,m,l\}$. Since $n\neq k$,
$m\neq l$, $k\neq l$ and $m\neq n$, we see that $j\geq 2$ and no
cancelation occurs in the expression $(r_n-r_k)-(r_m-r_l)$. Thus,
$$
|(r_n-r_k)-(r_m-r_l)|\geq r_j-r_{j-1}-r_{j-2}>a_j+a_{j-1}\geq
a_k+a_l.
$$

{\bf Case} $m=n$. Since $k>l$, we have
$$
|(r_n-r_k)-(r_m-r_l)|=r_k-r_l\geq r_k-r_{k-1}> 2a_n=a_n+a_m.
$$

Thus, (\ref{44}) is satisfied and therefore the supports of $y_k$
are disjoint. \square

Let $T$ be an invertible bilateral weighted shift acting on
$\ell_p(\Z)$, $1<p<\infty$ and $w=\{w_n\}_{n\in\Z}$ be its weight
sequence. As usual $\beta(a,b)$ stand for the numbers defined in
(\ref{beta}).

\prop{3.4}Suppose that there exist a sequence $\{r_n\}_{n\in\zp}$ of
positive integers and sequences $\{\alpha_n\}_{n\in\zp}$,
$\{\rho_n\}_{n\in\zp}$ of positive numbers such that
\begin{itemize}\itemsep=-2pt
\item[{\rm (W1)}]$\alpha_n\to\infty$ as $n\to\infty$;
\item[{\rm (W2)}]$r_{n+2}-r_{n+1}-r_n\to\infty$ as $n\to\infty$;
\item[{\rm (W3)}]$\sum\limits_{n=0}^\infty\rho^p_n\alpha^p_n
\beta(1,r_n)^{-p}<\infty$;
\item[{\rm (W4)}]$\sum\limits_{k=1}^\infty\biggl(\max\limits_{1\leq m\leq
k}\biggl(\sum\limits_{n=0}^{m-1}
\frac{\alpha^p_n\rho^p_n}{\rho^p_m}\beta(r_n-r_m+1,0)^p+
\sum\limits_{n=m+1}^\infty\frac{\alpha_n^p\rho^p_n}{\rho^p_m}
\beta(1,r_n-r_m)^{-p}\biggr)\biggr)^{-\frac1{p-1}}=\infty$.
\end{itemize}
Then $T$ is weakly supercyclic.

If a sequence $\{r_n\}_{n\in\zp}$ of positive integers and a
sequence $\{\alpha_n\}_{n\in\zp}$ of positive numbers can be chosen
such that conditions {\rm (W1--W4)} are satisfied with
$\rho_n\equiv1$, then $T$ is weakly hypercyclic. \rm

\proof Since $T$ is invertible, there exists $c>1$ such that
$c^{-1}\leq|w_n|\leq c$ for each $n\in\Z$. Therefore
$c^{a-b-1}\leq\beta(a,b)\leq c^{b-a+1}$ for each $a,b\in\Z$, $a\leq
b$. Moreover,
\begin{equation}
\frac{\beta(a,b)}{\beta(a+j,b+j)}\leq c^{2|j|}\ \ \text{for each
$a,b,j\in\Z$, $a\leq b$}. \label{wb2}
\end{equation}

Let $x$ be a vector from $\ell_p(\Z)$ with finite support. Using
(\ref{wb2}), we obtain that for any $n\in\N$,
\begin{align}
\|T^{-n}x\|_p&\leq \|x\|_p \max_{|j|\leq
\gamma(x)}(\beta(j+1,j+n))^{-1}\leq
\frac{\|x\|_p}{\beta(1,n)}\max_{|j|\leq
\gamma(x)}\frac{\beta(1,n)}{\beta(j+1,j+n)}\leq
\frac{\|x\|_p\,c^{2\gamma(x)}}{\beta(1,n)}; \label{minus}
\\
\|T^nx\|_p&\leq \|x\|_p \max_{|j|\leq \gamma(x)}\beta(j-n+1,j)\leq
\notag
\\
&\qquad\leq \|x\|_p\beta(1-n,0)\max_{|j|\leq
\gamma(x)}\frac{\beta(j-n+1,j)}{\beta(1-n,0)} \leq \|x\|_p
\beta(1-n,0)c^{2\gamma(x)}.\label{plus}
\end{align}

Note that the conditions (W1--W4) remain valid if we replace $r_n$,
$\alpha_n$ and $\rho_n$ by $r_{n+m}$, $\alpha_{n+m}$ and
$\rho_{n+m}$ respectively for any fixed non-negative integer $m$.
Thus, taking (W1) into account, we can, without loss of generality,
assume that $r_1>r_0$ and $r_{n}>r_{n-1}+r_{n-2}$ for each $n\geq
2$.

According to (W3) we have
$$
\sum_{n=k+1}^\infty\rho^p_n\alpha^p_n\beta(1,r_n-r_k)^{-p}\leq
c^{r_k}\sum_{n=k+1}^\infty\rho_n^p\alpha_n^p\beta(1,r_n)^{-p}<\infty
\ \ \text{for each}\ \ k\in\zp.
$$
Hence we can pick a strictly increasing sequence $\{m_k\}_{k\in\zp}$
of positive integers such that
\begin{equation}
\sum_{n=m_k}^\infty\rho^p_n\alpha^p_n\beta(1,r_n-r_k)^{-p}<\rho_k^p2^{-pk}
\ \ \text{for each}\ \ k\in\zp. \label{tail1}
\end{equation}

Now choose a monotonically non-decreasing sequence $\{a_n:n\in\zp\}$
of non-negative integers tending to infinity slowly enough to ensure
that
\begin{align}
2a_{m_k}&< r_k-r_{k-1}\ \ \text{for each $k\in\N$};\label{an1}
\\
a_nc^{2a_n}&\leq \alpha_n\ \ \text{for each $n\in\zp$}; \label{an2}
\\
a_n+a_{n-1}&<r_n-r_{n-1}-r_{n-2}\ \ \text{for each $n\geq 2$.}
\label{an3}
\end{align}

According to Lemma~3.2 there exists a map $\kappa:\zp\to\ell_p(\Z)$
such that the conditions (U1), (U2) and (U3) are satisfied. Since
$\supp(T^{-r_n}\kappa(n))\subseteq L(r_n,a_n)$, from (\ref{an3}) and
(\ref{ldis}) it follows that the supports of $T^{-r_n}\kappa(n)$ are
disjoint. The estimates (\ref{minus}), (U2) and (\ref{an2}) imply
that
$$
\|T^{-r_n}\kappa(n)\|_p\leq
a_nc^{2a_n}\beta(1,r_n)^{-1}\leq\alpha_n\beta(1,r_n)^{-1} \ \
\text{for each $n\in\zp$.}
$$
By (W3), $\sum\limits_{n=0}^\infty
\rho^p_n\|T^{-r_n}\kappa(n)\|_p^p<\infty$. Since the supports of
$T^{-r_n}\kappa(n)$ are disjoint, the series
$$
u=\sum_{n=0}^\infty \rho_nT^{-r_n}\kappa(n)
$$
is norm-convergent in $\ell_p(\Z)$. It suffices to prove that $u$ is
a weakly supercyclic vector for $T$ and that $u$ is a weakly
hypercyclic vector for $T$ if $\rho_n\equiv1$.

Clearly $T^{r_k}u=\rho_k\kappa(k)+v_k+z_k+y_k$, where
$$
v_k=\sum_{a_n\geq (r_k-r_{k-1})/2\atop n>k}\!\!
\rho_nT^{r_k-r_n}\kappa(n),\ \ z_k=\sum_{a_n<(r_k-r_{k-1})/2\atop
n>k}\!\! \rho_nT^{r_k-r_n}\kappa(n), \ \ y_k=\sum_{n=0}^{k-1}
\rho_nT^{r_k-r_n}\kappa(n).
$$
From (\ref{plus}), (\ref{minus}), (U2) and (\ref{an2}), we have
\begin{alignat}{2}
\|T^{r_k-r_n}\kappa(n)\|_p&\leq a_nc^{2a_n}\beta(1,r_n-r_k)^{-1}\leq
\alpha_n\beta(1,r_n-r_k)^{-1},& \quad &\text{if $n>k$}; \label{nbi}
\\
\|T^{r_k-r_n}\kappa(n)\|_p&\leq a_nc^{2a_n}\beta(r_n-r_k+1,0)\leq
\alpha_n\beta(r_n-r_k+1,0),& \quad &\text{if $n<k$}. \label{nsm}
\end{alignat}

Since $T$ preserves disjointness of the supports and the supports of
$T^{-r_n}\kappa(n)$ are disjoint, we see that for any $k\in\zp$ the
supports of $T^{r_k-r_n}\kappa(n)$, $n\in\zp$ are also disjoint.
Hence,
$$
\|v_k\|_p^p=\sum_{a_n\geq (r_k-r_{k-1})/2\atop n>k}\!\!\!\!\!
\rho_n^p\|T^{r_k-r_n}\kappa(n)\|_p^p\ \ \text{for any $k\in\N$}.
$$
Applying (\ref{an1}), (\ref{tail1}) and (\ref{nbi}), we obtain
\begin{equation}
\|v_k\|_p^p\leq \sum_{n\geq
m_k}\rho_n^p\|T^{r_k-r_n}\kappa(n)\|_p^p\leq \sum_{n\geq
m_k}\rho^p_n \alpha^p_n \beta(1,r_n-r_k)^{-p}\leq \rho^p_k2^{-kp}\ \
\text{for any $k\in\N$.} \label{vk}
\end{equation}
Analogously, applying (\ref{nbi}) and (\ref{nsm}), we get
\begin{align*}
\|z_k\|^p_p&\leq \sum_{n=k+1}^\infty
\rho_n^p\|T^{r_k-r_n}\kappa(n)\|^p_p\leq \sum_{n=k+1}^\infty
\alpha_n^p\rho_n^p\beta(1,r_n-r_k)^{-p};
\\
\|y_k\|_p^p&=\sum_{n=0}^{k-1}\rho^p_n\|T^{r_k-r_n}\kappa(n)\|^p_p\leq
\sum_{n=0}^{k-1} \alpha^p_n\rho^p_n\beta(r_n-r_k+1,0)^{p}.
\end{align*}
Hence
\begin{align*}
\|z_k+y_k\|^p_p&=\|y_k\|_p^p+\|z_k\|_p^p\leq \rho^p_k\xi_k\leq
\rho^p_k\theta_k,\ \ \ \text{where}
\\
\xi_m&=\sum_{n=0}^{m-1}
\frac{\alpha_n^p\rho_n^p}{\rho_m^p}\beta(r_n-r_m+1,0)^{p}+
\sum_{n=m+1}^\infty
\frac{\alpha_n^p\rho_n^p}{\rho_m^p}\beta(1,r_n-r_m)^{-p}\ \
\text{and}\ \ \theta_k=\max_{1\leq m\leq k} \xi_m.
\end{align*}

In view of (\ref{an3}), Lemma~3.3 implies that the sequence
$\{y_k+z_k\}_{k\in\N}$ has disjoint supports. Hence
$\{\rho_k^{-1}\theta_k^{-1/p}(y_k+z_k)\}$ is a $p$-sequence in
$\ell_p(\Z)$ as a bounded sequence with disjoint supports. Since the
sequence $\theta_k$ is monotonically non-decreasing, it is bounded
from below by a positive constant and therefore from (\ref{vk}) and
Lemma~2.5 on small perturbations of $p$-sequences it follows that
$\{\rho_k^{-1}\theta_k^{-1/p}(v_k+y_k+z_k)\}$ is a $p$-sequence in
$\ell_p(\Z)$. Since $T^{r_k}u=\rho_k\kappa(k)+v_k+z_k+y_k$, we see
that
$\{\rho_k^{-1}\theta_k^{-1/p}T^{r_k}u-\theta_k^{-1/p}\kappa(k)\}_{k\in\N}$
is a $p$-sequence in $\ell_p(\Z)$.

Pick $x\in S=\kappa(\zp)\setminus\{0\}$ and let
$A_x=\kappa^{-1}(x)$. By (U3) $A_x$ has positive density. According
to (W4), $\sum\limits_{k=1}^\infty \theta_k^{-1/(p-1)}=\infty$. By
Lemma~3.1 $\sum\limits_{k\in A_x} \theta_k^{-1/(p-1)}=\infty$, or
equivalently, $\sum\limits_{k\in A_x} (\theta_k^{-1/p})^q=\infty$,
where $\frac1p+\frac1q=1$. Since $\kappa(k)=x$ for any $k\in A_x$,
we see that
$\{\rho_k^{-1}\theta_k^{-1/p}T^{r_k}u-\theta_k^{-1/p}x\}_{k\in A_x}$
is a $p$-sequence in $\ell_p(\Z)$. Since $S$ is dense in
$\ell_p(\Z)$, Proposition~2.4 implies that $u$ is a weakly
supercyclic vector for $T$. If additionally $\rho_n\equiv1$,
Proposition~2.4 implies that $u$ is a weakly hypercyclic vector for
$T$. \square

{\bf Remark.} \ An analog of Proposition~3.4 holds for invertible
bilateral weighted shifts on $c_0(\Z)$ and the proof is basically
the same. One has to replace conditions (W3) and (W4) by
\begin{itemize}\itemsep=-2pt
\item[{\rm (W$3'$)}]$\lim\limits_{n\to\infty}\rho_n\alpha_n
\beta(1,r_n)^{-1}=0$;
\item[{\rm (W$4'$)}]$\sum\limits_{k=1}^\infty \biggl(
\max\limits_{1\leq m\leq k}\biggl(\max\limits_{0\leq n\leq m-1}
\frac{\alpha_n\rho_n}{\rho_m}\beta(r_n-r_m+1,0)+ \max\limits_{n\geq
m+1}\frac{\alpha_n\rho_n}{\rho_m}
\beta(1,r_n-r_m)^{-1}\biggr)\biggr)^{-1}=\infty$.
\end{itemize}

\subsection{Proof of Theorem~1.5}

We have to prove that $B$ is weakly supercyclic for $p>2$. Take
$\{r_k\}$ being any sequence of positive integers, satisfying the
condition (W2) of Proposition~3.4, for instance $r_k=2^k$. Clearly
for the unweighted bilateral shift, we have $\beta(a,b)=1$ for each
$a,b\in\Z$, $a\leq b$. Therefore, if (W3) is satisfied, then the
$k$-th term in the sum in (W4) is bounded from below by
$c\rho_k^{p/(p-1)}$ for some positive constant $c$. Thus, all
conditions of Proposition~3.4 will be satisfied if we find sequences
$\{\alpha_n\}_{n\in\zp}$ and $\{\rho_n\}_{n\in\zp}$ of positive
numbers such that $\alpha_n\to\infty$, $\sum\limits_{n=0}^\infty
\alpha^p_n\rho^p_n<\infty$ and $\sum\limits_{k=0}^\infty
\rho_k^{p/(p-1)}=\infty$. This can be achieved by choosing
$\alpha_n=\ln(n+2)$ and $\rho_n=(n+1)^{-1/p}(\ln(n+2))^{-2}$.

\subsection{Proof of Proposition 1.8}

Consider the sequence $w=\{w_n\}_{n\in\Z}$ defined by the formula
$$
w_m=\left\{\begin{array}{ll}2&\text{if $7\cdot9^k<m\leq9^{k+1}$, $k$
even, or $-11\cdot9^k\leq m<-9^{k+1}$, $k$ odd};
\\
1/2&\text{if $9^{k+1}\leq m<11\cdot9^{k+1}$, $k$ even, or
$-9^{k+1}<m\leq-7\cdot9^k$, $k$ odd};
\\
1&\text{otherwise.}\end{array} \right.
$$
In this section $T$ stands for the bilateral weighted shift with the
weight sequence $w$, acting on $\ell_2(\Z)$. Obviously $T$ is
invertible. From definition of the weight sequence $\{w_n\}$ it
follows that $\max\{\beta(-n,0),(\beta(0,n))^{-1}\}\geq1$ for each
$n\in\zp$. Hence $T$ is not hypercyclic according to Theorem~S. It
remains to show that $T$ and $T^{-1}$ are weakly hypercyclic.

Consider the sequences $\rho_n=1$, $r_n=9^{2n+1}$ and
$\alpha(n)=\ln\ln(n+4)$, $n\in\zp$. Conditions (W1) and (W2) of
Proposition~3.4 are trivially satisfied.  Using the definition of
the weight sequence $\{w_n\}$, one can easily verify that for even
$k\in\zp$
\begin{equation}
\beta(1,a)=2^{a-7\cdot9^k}\ \ \text{and}\ \ \beta(-a+1,0)=1\ \
\text{if $7\cdot9^k<a\leq 9^{k+1}$}. \label{bebeta}
\end{equation}
This implies that $\beta(1,r_n-r_k)>2^{9^{2n}}$ if $n>k$ and
$\beta(r_n-r_k+1,0)=1$ if $n<k$. Now it is an elementary exercise to
show that (W3) and (W4) for $p=2$ are also satisfied. Thus, $T$ is
weakly hypercyclic according to Proposition~3.4.

Since the operator $T^{-1}$ is similar to the weighted bilateral
shift $\widetilde T$ with the weight sequence $\widetilde
w_n=w_{-n}^{-1}$, it suffices to verify that $\widetilde T$ is
weakly hypercyclic. This follows from Proposition~3.4 similarly via
choosing the sequences $\rho_n=1$, $r_n=9^{2n+2}$ and
$\alpha(n)=\ln\ln(n+4)$, $n\in\zp$. Indeed, (\ref{bebeta}) is
satisfied for odd $k\in\zp$ for the weight sequence $\widetilde
w_n$.

\subsection{Proof of Proposition 1.9}

Let $p>2$ and $\phi:[0,\infty)\to [1,\infty)$ be the function
defined as
$$
\phi(t)=(t+1)^{1/p}(\log_2(t+2))^{2/p}.
$$
Consider the sequence $w=\{w_n\}_{n\in\Z}$ of positive numbers
defined by the formula
$$
w_m=\left\{\begin{array}{ll}\left(\frac{\phi(k+1)}{\phi(k)}
\right)^{3^{-k}/2}&\text{if}\ 3^n-3^{k+1}<|m|\leq 3^n-3^k,\
n,k\in\zp,\ n\geq k+2; \\
\left(\frac{\phi(n+1)}{\phi(n)^2} \right)^{3^{-n}}&\text{if}\
3^n<|m|\leq 3^{n+1}-3^n,\
n\in\zp; \\
1&\text{if $m=0$ or $|m|=3^n$, $n\in\zp$}.\end{array}\right.
$$
In this section $T$ stands for the bilateral weighted shift with the
weight sequence $w$, acting on $\ell_p(\Z)$. It suffices to prove
that $T$ is weakly hypercyclic and non-supercyclic.

It is easy to see that the sequence $w$ is symmetric: $w_n=w_{-n}$,
$n\in\zp$ and that $w_n\to1$ as $|n|\to\infty$. Hence there exists
$c>1$ such that $c^{-1}\leq w_n\leq c$ for each $n\in\Z$. Therefore
$T$ is invertible. Using the definition of $w$ it is straightforward
to verify that
\begin{align}
\beta(1,3^n-3^k)&=\phi(n)/\phi(k)\ \text{if $n,k\in\zp$ and $n>k$},
\label{be1}
\\
\beta(1,3^n)&=\phi(n)\ \text{for $n\in\zp$}. \label{be2}
\end{align}
We shall prove that $T$ is weakly hypercyclic. Consider the
sequences $\rho_n=1$, $r_n=3^n$ and $\alpha(n)=\ln\ln(n+4)$,
$n\in\zp$. Conditions (W1) and (W2) of Proposition~3.4 are trivially
satisfied. Using (\ref{be2}), we see that
$$
\sum_{n=0}^\infty
\rho^p_n\alpha^p_n\beta(1,r_n)^{-p}=\sum_{n=0}^\infty\alpha^p_n\phi(n)^{-p}
<\infty.
$$
Hence (W3) is also satisfied. Let now
$$
\xi_m=\sum_{n=0}^{m-1}\frac{\alpha^p_n\rho^p_n}{\rho^p_m}\beta(r_n-r_m+1,0)^p+
\sum_{n=m+1}^\infty\frac{\alpha^p_n\rho^p_n}{\rho^p_m}\beta(1,r_n-r_m)^{-p}\
\ \text{for $m\in\N$}.
$$
Since $\rho_n\equiv1$ and the weight sequence $w$ is symmetric, we
using (\ref{be1}) obtain
\begin{align*}
\xi_m&=\sum_{n=0}^{m-1}\alpha^p_n\beta(1,r_m-r_n)^p\frac{w_0^p}{w_{r_n-r_m}^p}+
\sum_{n=m+1}^\infty \alpha^p_n\beta(1,r_n-r_m)^{-p}\leq
\\
&\leq c^p\sum_{n=0}^{m-1}\alpha^p_n\frac{\phi(m)^p}{\phi(n)^p}+
\sum_{n=m+1}^\infty\alpha^p_n\frac{\phi(m)^p}{\phi(n)^p}\leq
c^p\phi(m)^p\sum_{n=0}^\infty\alpha_n^p\phi(n)^{-p}=A(p)\phi(m)^p,
\end{align*}
where $A(p)$ is a positive constant depending only on $p$. Since
$p>2$, we have
$$
\sum_{k=1}^\infty\Bigl(\max\limits_{1\leq m\leq
k}\xi_m\Bigr)^{-\frac1{p-1}}\geq A(p)^{-1/(p-1)}\sum_{k=1}^\infty
\phi(k)^{-p/(p-1)}=\infty.
$$
Thus, (W4) is also satisfied and Proposition~3.4 implies that $T$ is
weakly hypercyclic. It remains to notice that according to
Theorem~S, a bilateral weighted shift with symmetric weight sequence
is never supercyclic. The proof is complete.

\section{The multiplication operator $M$: proof of Theorem~1.2}

Let $\M=\M(\T)$ be the space of $\sigma$-additive complex-valued
Borel measures on the unit circle $\T$. We denote the set of
non-negative measures $\mu\in\M$ as $\M_+$. It is well-known that
$\M$ is a Banach space with respect to the variation norm
$\|\mu\|=|\mu|(\T)$, where $|\mu|\in\M_+$ is the variation of $\mu$.
That is, $|\mu|(A)$ is the supremum of $\sum\limits_n|\mu(A_n)|$,
where $A_n$ are disjoint Borel subsets of $A$. The set of measures
$\mu\in\M$, whose Fourier coefficients $
\widehat\mu(n)=\int\limits_{\T}z^n\,d\mu(z)$, $n\in\Z$ tend to zero
when $|n|\to\infty$ will be denoted by $\M_0$.

\subsection{Proof of Proposition~1.1}

Let $f,g\in L_2(\mu)$. Then $\langle M^nf,g\rangle=\widehat\nu(n)$,
where $\nu\in\M$ is absolutely continuous with respect  to $\mu$
with the density $\frac{d\nu}{d\mu}(z)=f(z)\overline{g(z)}$. Since a
measure absolutely continuous with respect to a measure from $\M_0$
also belongs to $\M_0$, see \cite{gg}, we find that $\nu\in\M_0$,
that is, $\widehat\nu(n)\to0$. Hence $\langle
M^nf,g\rangle/\|M^nf\|=\widehat\nu(n)/\|f\|\to 0$ for any non-zero
$f\in L_2(\mu)$ and any $g\in L_2(\mu)$. Therefore the sequence
$\{M^nf/\|M^nf\|\}$ tends weakly to zero. Thus, $M$ is
antisupercyclic. It remains to apply Theorem~2.1.

\subsection{Weak convergence of measures}

We need to introduce further notation. For a Borel measurable set
$K\subset\T$ we denote by $\M(K)$ the set of $\mu\in\M$ such that
$|\mu|(\T\setminus K)=0$. The {\it support} of a measure $\mu\in\M$
is
$$
\supp(\mu)=\bigcap\{K\subset\T:K\ \text{is closed and}\ \mu\in\M(K)\}.
$$
Recall that the weak topology $\sigma$ on $\M$  is the topology
generated by the functionals
$$
\mu\mapsto[\mu,f]=\int\limits_\T f(z)\,\mu(dz), \quad f\in {\cal
C}(\T),
$$
that is, $\sigma$ is the weakest topology with respect  to which the
functionals $\mu\mapsto[\mu,f]$ are continuous. For $\mu\in\M$ and a
Borel-measurable set $A\subset \T$, $\ssub\mu{\!A}$ stands for the
restriction of $\mu$ to $A$, that is $\ssub\mu{\!A}\in\M$ is defined
by $\ssub\mu{\!A}(B)=\mu(A\cap B)$. An {\it interval} $I$ of $\T$ is
a non-empty {\bf open}  connected subset of $\T$ and $|I|$ will
denote its length.

\lemma{4.1} Let $\I^n=\{I_1^n,I_2^n,\dots,I_{k_n}^n\}$, $n\in\N$ be
a family of disjoint intervals of $\T$ satisfying
\begin{itemize}\itemsep=-2pt
\item[{\rm (L1)}]for each $n\in\N$, any element of $\I^{n+1}$ is
contained in some element of $\I^n$;
\item[{\rm (L2)}]$\max\limits_{1\leq j\leq k_n}|I^n_j|\to 0$ as
$n\to\infty$.
\end{itemize}
Let also $\mu\in\M_+$ and $\mu^n\in\M_+$, $n\in\N$ be such that
\begin{itemize}\itemsep=-2pt
\item[{\rm (L3)}]$\mu(\T)=\mu^n(\T)=
\sum\limits_{j=1}^{k_n}\mu^m(I^n_j)=\sum\limits_{j=1}^{k_n}\mu(I^n_j)$
for any $n\in\N$; \item[{\rm (L4)}]$\mu(I^n_j)=\mu^n(I^n_j)$ for any
$n\in\N$ and $j=1,\dots,k_n$.
\end{itemize}
Then the $\mu^n\tos\mu$ as $n\to\infty$. Moreover
$\ssub{\mu^n}{I^m_j}\tos \ssub{\mu}{I^m_j}$ as $n\to\infty$ for any
$m\in\N$ and $1\leq j\leq k_m$. \epr

\proof Let $f\in{\cal C}(\T)$ and $\epsilon>0$. Since $f$ is
uniformly continuous, condition (L2) implies the existence of
$a\in\N$ and $z_1,\dots, z_{k_a}\in\C$ for which $|f(z)-z_j|\leq
\epsilon$ if $z\in I^a_j$ and $1\leq j\leq k_a$. According to (L1),
(L3) and (L4), $\mu^n(I_j^a)=\mu(I_j^a)\ \ \text{for each $n\geq a$
and $1\leq j\leq k_a$}$.  Hence for each $n\geq a$, we have
$$
|[\ssub{\mu}{I^a_j},f]-c_jz_j|\leq c_j\epsilon\ \ \text{and}\ \
|[\ssub{\mu^n}{I^a_j},f]-c_jz_j|\leq c_j\epsilon\text{\ \ for $1\le
j\leq n_a$, where $c_j=\mu(I_j^a)$.}
$$
Therefore $ |[\ssub{\mu}{I^a_j},f]-[\ssub{\mu^n}{I^a_j},f]|\leq
2c_j\epsilon\ \ \text{if $n\geq a$ and $1\leq j\leq k_a$}$. Thus
summing over $j$, we obtain
$$
|[\mu,f]-[\mu^n,f]|\leq
2\epsilon\sum_{j=1}^{k_a}c_j=2\mu(\T)\epsilon \ \ \text{for}\  \
n\geq a.
$$
Hence $\mu^n\tos\mu$ as $n\to\infty$.

Fix now $m,j\in\N$ such that $1\leq j\leq k_m$. One can easily
verify that conditions (L1--L4) remain valid if we replace $\I^n$ by
$\I^{n+m}$, $\mu^n$ by $\ssub{\mu^{n+m}}{I^m_j}$ and $\mu$ by
$\ssub{\mu}{I^m_j}$. From what is already proven, follows that
$\ssub{\mu^{n+m}}{I^m_j}\tos \ssub{\mu}{I^m_j}$. Hence
$\ssub{\mu^n}{I^m_j}\tos \ssub{\mu}{I^m_j}$. \square

Recall that a set $A\subset\T$ is called {\it independent} if
$\prod\limits_{j=1}^m z_j^{n_j}\neq 1$ for each pairwise different
points $z_1,\dots,z_m\in A$ and each non-zero vector
$(n_1,\dots,n_m)\in\Z^m$.

The set of probability measures $\mu\in\M$, will be denoted by $\P$,
$\pac$ will denote the set of measures in $\P$ absolutely continuous
with respect to the Lebesgue measure and $\pfin$ will denote the set
of measures in $\P$ with finite independent support.

\lemma{4.2} Let  $I_1,\dots,I_m$  be disjoint intervals of $\T$,
$A=\bigcup\limits_{j=1}^m I_j$ and $\mu\in\P\cap\M(A)$. Then there
exist sequences $\mu^n\in\pac$ and $\nu^n\in\pfin$ $(n\in\N)$ such
that
\begin{align}
&\mu^n(I_j)=\nu^n(I_j)=\mu(I_j)\ \ \text{for}\  \ 1\leq j\leq m\ \
\text{and} \label{du1}
\\
&\mu^n\tos\mu,\  \nu^n\tos\nu,\
\ssub{\mu^n}{I_j}\tos\ssub{\mu}{I_j},\
\ssub{\nu^n}{I_j}\tos\ssub{\mu}{I_j}\ \text{as $n\to\infty$ for
$1\leq j\leq m$}. \label{du2}
\end{align}
\rm

\proof Let $B$ be the set of atoms of $\mu$, that is
$B=\{z\in\T:\mu(\{z\})>0\}$, which is at  most countable since $\mu$
is finite. Thus, for each interval $J$ of $T$ and $\epsilon>0$,
there exists a disjoint family of intervals $J_1,\dots,J_d$ such
that $\max\limits_{1\leq j\leq d}|J_j|<\epsilon$, $J_j\subset J$ for
any $j=1,\dots, d$ and $J\setminus\bigcup \limits_{j=1}^d J_j$ is
finite and does not meet $B$. In this way it is  easy to choose a
sequence $\I^n=\{I_1^n,I_2^n,\dots,I_{k_n}^m\}$ of disjoint families
of intervals of $\T$ with $\I^1=\{I_1,\dots,I_m\}$, for which (L1)
and (L2) are satisfied and for any $n\in\N$ and the set
$\bigcup\limits_{j=1}^{k_n}
I_j^{n}\setminus\bigcup\limits_{l=1}^{k_{n+1}}I_l^{n+1}$ is finite
and does not meet $B$. The latter property implies that
\begin{equation}
\sum\limits_{j=1}^{k_n}\mu(I^n_j)=1\ \ \text{for each}\ \ n\in\N.
\label{qqq}
\end{equation}
Now we can define $\mu^n$ and $\nu^n$. Let $\mu^n$ be the absolutely
continuous measure with the density
$$
\rho^n(z)=\sum_{j=1}^{k_n} c^n_j\ssub{\chi}{I^n_j},
$$
where $\ssub\chi I$  denotes the indicator function of a set $I$ and
$c^n_j\geq0$ are chosen in such a way that $\mu^n(I^n_j)=\mu(I^n_j)$
for $1\leq j\leq k_n$. Hence $\mu^n(I_j)=\mu(I_j)$ for $n\in\N$ and
$1\leq j\leq m$. From (\ref{qqq}) it also follows that
$\mu^n\in\pac$.

To define $\nu^n$, choose a set of independent points $z_{n,j}\in
I^n_j$, $1\leq j\leq k_n$ and consider
$$
\nu^n=\sum_{j=1}^{k_n}\mu(I^n_j)\delta_{z_{n,j}},
$$
where $\delta_z$ stands for the probability measure with the
one-point support $\{z\}$. From (\ref{qqq}), we see that
$\nu^n\in\P$ and therefore $\nu^n\in\pfin$. Obviously
$\nu^n(I^n_j)=\mu(I^n_j)$ for $1\leq j\leq k_n$ and therefore
$\nu^n(I_j)=\mu(I_j)$ for $j=1,\dots,m$. Thus, (\ref{du1}) holds.
Finally (\ref{du2})  follows  from Lemma~4.1. \square

Next lemma is the main  building block in the inductive procedure of
constructing the measure, asserted by Theorem~1.2.

\lemma{4.3}Let $\epsilon>0$, $k_0\in\N$, $h_1,\dots,h_n\in {\cal
C}(\T)$, $I_1,\dots,I_m$ be disjoint intervals,
$c_1,\dots,c_m\in\C\setminus\{0\}$ with
$a=\max\{|c_1|,\dots,|c_m|\}\leq 1$ and $\mu\in\pac$ be such that
$\sum\limits_{j=1}^m \mu(I_j)=1$. Then there exist $\nu\in\pac$ and
$k\in\N$, $k>k_0$ satisfying
\begin{itemize}\itemsep=-2pt
\item[{\rm(B1)}] $\nu(I_j)=\mu(I_j)$ \ for $1\leq j\leq m$;
\item[{\rm(B2)}] $|[\mu-\nu,h_l]|<\epsilon$ for $1\leq l\leq n$;
\item[{\rm(B3)}] $|\widehat{\ssub{\nu}{I_j}}(k)-
c_j\mu(I_j)|<\epsilon$ for $1\leq j\leq m$;
\item[{\rm(B4)}] $\|\mmu-\nnu\|_\infty\leq 2a$.
\end{itemize}
\epr

\proof Since $\ssub{\mu}{I_j}$ are absolutely continuous,
$\widehat{\ssub{\mu}{I_j}}(k)\to 0$ as $|k|\to\infty$ for $1\leq
j\leq m$. Therefore there exists $k_1\in\N$ such that
$|\widehat{\ssub{\mu}{I_j}}(k)|<\epsilon/3$ if $|k|\geq k_1$ and
$1\leq j\leq m$. By Lemma~4.2, there exists $\gamma\in\pfin$ such
that
\begin{align}
\gamma(I_j)=\mu(I_j)\ \ &\text{for}\ \ 1\leq j\leq m; \label{435}
\\
|[\ssub{\mu}{I_j}-\ssub{\gamma}{I_j},h_l]|< \epsilon/(2m)\ \
&\text{for}\ \ 1\leq j\leq m\ \text{and}\ 1\leq l\leq n. \label{436}
\end{align}
Since $\supp(\gamma)=\{u_1,\dots,u_N\}$ is an independent set, the
Kronecker theorem \cite{gg} implies that the sequence
$\{(u_1^k,\dots,u_N^k)\}:k\in\N\}$ is dense in $\T^N$. Consider the
vector $\theta\in\T^N$ defined by
$$
\theta_s=c_j/|c_j|\ \ \text{if}\  \ u_s\in I_j.
$$
Choosing $k>\max\{k_0,k_1\}$ in such a way that
$(u_1^k,\dots,u_N^k)$ is close enough to $\theta$, we can ensure
that
\begin{equation}
|\widehat{\ssub{\gamma}{I_j}}(k)-c_j\gamma(I_j)/|c_j||=
|\widehat{\ssub{\gamma}{I_j}}(k)-c_j\mu(I_j)/|c_j||<\epsilon/3\ \
\text{for}\ \ 1\leq j\leq m. \label{02}
\end{equation}
Applying Lemma~4.2 once again, we
obtain that there exists $\eta\in\pac$ such that
\begin{align}
\eta(I_j)=\gamma(I_j)=\mu(I_j)\ \ &\text{for}\ \ 1\leq j\leq
m;\label{437}
\\
|[\ssub{\eta}{I_j}-\ssub{\gamma}{I_j},g_l]|< \epsilon/(2m)\ \
&\text{for \ $1\leq j\leq m$ and $1\leq l\leq n$}; \label{438}
\\
|\widehat{\ssub{\eta}{I_j}}(k)-
\widehat{\ssub{\gamma}{I_j}}(k)|<\epsilon/3\ \ &\text{for}\ \ 1\leq
j\leq m. \label{439}
\end{align}
The required measure is
\begin{equation}
\nu=\sum_{j=1}^m(1-|c_j|)\ssub{\mu}{I_j}+|c_j|\ssub{\eta}{I_j},
\label{defnu}
\end{equation}
which, since $|c_j|\leq 1$, is non-negative and clearly absolutely
continuous. From (\ref{437}) we find that $\nu(I_j)=\mu(I_j)$ for
$1\leq j\leq m$. Thus, $\nu\in\pac$ and (B1) holds. From
(\ref{defnu}) it follows that
$$
|[\mu-\nu,g_l]|=\left|\sum_{j=1}^m|c_j|[\ssub{\mu}{I_j}-
\ssub{\eta}{I_j},g_l]\right|\leq a\sum_{j=1}^m(|[\ssub{\eta}{I_j}-
\ssub{\gamma}{I_j},g_l]|+
|[\ssub{\mu}{I_j}-\ssub{\gamma}{I_j},g_l]|).
$$
Using (\ref{436}) and (\ref{438}), we obtain
$$
|[\mu-\nu,g_l]|< am(\epsilon/(2m)+\epsilon/(2m))\leq a\epsilon\leq
\epsilon\ \ \text{for $1\leq j\leq n$}.
$$
Hence (B2) holds. Suppose now that $1\leq j\leq m$. Then
\begin{align*}
|\widehat{\ssub{\nu}{I_j}}(k)-c_j\mu(I_j)|&=|(1-|c_j|)
\widehat{\ssub{\mu}{I_j}}(k)+
|c_j|\widehat{\ssub{\eta}{I_j}}(k)-
c_j\mu(I_j)|\leq \\ &\leq|\widehat{\ssub{\mu}{I_j}}(k)|+
|c_j||\widehat{\ssub{\eta}{I_j}}(k)-
\widehat{\ssub{\gamma}{I_j}}(k)|+
|c_j||\widehat{\ssub{\gamma}{I_j}}(k)-c_j\mu(I_j)/|c_j||.
\end{align*}
Thus by (\ref{439}), (\ref{02}) and the inequalities
$|\widehat{\ssub{\mu}{I_j}}(k)|<\epsilon/3$, $|c_j|\leq 1$, we
obtain $|\widehat{\ssub{\nu}{I_j}}(k)-c_j\mu(I_j)|<\epsilon$, that
is, (B3) holds. Finally from (\ref{437}) it follows that
$\|\ssub{\mu}{I_j}\|=\|\ssub{\eta}{I_j}\|=\mu(I_j)$ and we have
$$
\|\mmu-\nnu\|_\infty\leq
\|\mu-\nu\|=\left\|\sum_{j=1}^m |c_j|(\ssub{\mu}{I_j}-
\ssub{\eta}{I_j})\right\|\leq 2a
\sum_{j=1}^m \mu(I_j)=2a.
$$
Thus, (B4) also holds. \square

\lemma{4.4}Let $\delta_n>0$ and $f_n\in {\cal C}(\T)$ $(n\in\N)$ be
such that $\|f_n\|_\infty\to0$ as $n\to\infty$ and
$\|f_n\|_\infty\leq1$ for any $n\in\N$. Then there exists
$\mu\in\P\cap\M_0$ and a strictly increasing sequence $k_n$ of
non-negative integers such that
\begin{equation}
|[\mu,g_n\overline{g_m}]|\leq \delta_n\ \ \text{whenever}\ \ n>m,
\end{equation}
where $g_n(z)=z^{k_n}-f_n(z)$.
\epr

\proof First of all, we take a sequence $\{\epsilon_n\}_{n\in\N}$ of
positive numbers such that
\begin{equation}
\sum_{k=n}^\infty\epsilon_k<\delta_n/6\  \ \text{for each}\  \
n\in\N. \label{epde}
\end{equation}
For any $n,k\in\N$, $1\leq k\leq n$, let $I^n_k$ be the interval of
$\T$ between $e^{2\pi i(k-1)/n}$ and $e^{2\pi i k/n}$ (going
counterclockwise). Obviously, for any fixed $n$, the  intervals
$I^n_1,\dots,I^n_n$ are disjoint and
$\T\setminus\bigcup\limits_{j=1}^n I^n_j$ is finite. Therefore
$\nu=\sum\limits_{k=1}^n\ssub{\nu}{I^n_k}$ for any $n\in\N$ and any
continuous measure $\nu\in\M$.

Set $a_n=\|f_n\|_\infty$. For each $n\in\N$ we shall construct
inductively non-negative integers $k_n$, $j_n$, $m_n$, complex
numbers $c^{n,d}_j,b^n_j$ ($1\leq d\leq n$, $1\leq j\leq m_n$) and a
measure $\mu^n\in\pac$,  satisfying the following conditions:
\begin{itemize}\itemsep=-2pt
\item[(P1)]$0\leq k_{n-1}<k_n$, $1\leq j_{n-1}<j_n$ and $1\leq
m_{n-1}<m_n$ if $n\geq 2$; \item[(P2)]$m_{n-1}$ is a divisor of
$m_n$ for $n\geq 2$; \item[(P3)]$|c^{n,d}_j|\leq 2$, $|b^n_j|\leq
a_{n+1}$, $|g_d(z)-c^{n,d}_j|\leq\epsilon_{n+1}$ and
$|f_{n+1}(z)-b^n_j|\leq\epsilon_{n+1}$ for $1\leq d\leq n$, $1\leq
j\leq m_n$ and $z\in I^{m_n}_j$, where $g_d(z)=z^{k_d}-f_d(z)$;
\item[(P4)]$\mu^{n}(I^{m_{n-1}}_j)=\mu^{n-1}(I^{m_{n-1}}_j)$
for $1\leq j\leq m_{n-1}$ and $n\geq 2$;
\item[(P5)]$|\widehat{\mu^n_j}
(k_n)-b_j^{n-1}\mu^n(I^{m_{n-1}}_j)|\leq \epsilon_n/m_{n-1}$ for
$1\leq j\leq m_{n-1}$ and $n\geq 2$, where $\mu^n_j$ is the
restriction of $\mu^n$ to $I^{m_{n-1}}_j$;
\item[(P6)]$|\widehat{\mu^n}(l)|\leq \epsilon_n$ for $|l|\geq
j_n$;
\item[(P7)]$|\widehat{\mu^{n}}(l)-\widehat{\mu^{n-1}}(l)|\leq
\epsilon_n$ for $|l|<j_{n-1}$ if $n\geq2$;
\item[(P8)]$\|\widehat{\mu^{n}}-\widehat{\mu^{n-1}}\|_\infty
\leq 2a_{n+1}$ for $n\geq 2$;
\item[(P9)]$|[\mu^{n}-\mu^{n-1},g_m\overline{g_l}]|<\epsilon_n$
for $1\leq l<m\leq n-1$ and $n\geq 2$.
\end{itemize}

Set $k_1=0$ and $g_1(z)=z^{k_1}-f_1(z)=1-f_1(z)$. Take an arbitrary
measure  $\mu^1\in\pac$. Since the functions $f_2$ and $g_1$ are
uniformly continuous, there exist $m_1\in\N$ and complex numbers
$c^{1,1}_j$, $b^1_j$ ($1\leq j\leq m_1$) such that
$|c^{1,1}_j|\leq\|g_1\|_\infty\leq1+a_1\leq 2$,
$|b^n_j|\leq\|f_2\|_\infty=a_{2}$ and
$|g_1(z)-c^{1,1}_j|\leq\epsilon_2$, $|f_{2}(z)-b^1_j|\leq\epsilon_2$
if $1\leq j\leq m_1$ and $z\in I^{m_1}_j$. Since $\mu^1$ is
absolutely continuous, there exists $j_1\in\N$ such that
$|\widehat{\mu^1}(l)|\leq \epsilon_1$ for $|l|\geq j_1$. Thus,
$k_1$, $j_1$, $m_1$, $c^{1,1}_j,b^1_j$  and $\mu^1$ satisfy (P3) and
(P6): the only conditions required for $n=1$. The first step (basis)
of induction is done.

Assume now that $n\geq 2$ and $k_l$, $j_l$, $m_l$, $c^{l,d}_j$,
$b^l_j$ and $\mu^l$ for $1\leq l\leq n-1$, $1\leq d\leq l$ and
$1\leq j\leq m_l$, satisfying (P1--P9) are already constructed. We
have to construct $k_n$, $j_n$, $m_n$, $c^{n,d}_j$, $b^n_j$  and
$\mu^n$.

By Lemma~4.3 applied for $\mu=\mu^{n-1}$, $I_j=I^{m_{n-1}}_j$,
$c_j=b^{n-1}_j$, $k_0=k_{n-1}$, $m=m_{n-1}$,
$\epsilon=\epsilon_n/m_{n-1}$ and the finite set of functions
$\{h_j\}$ being $\{g_m\overline{g_l}:1\leq l<m\leq n-1\}\cup
\{z^l:|l|<j_{n-1}\}$, there exists a measure $\mu^n\in\pac$ and
$k_n\in\N$ such that $k_n>k_{n-1}$ and (P4), (P5), (P7), (P8) and
(P9) are satisfied.

Since $f_{n+1}$ and $g_d(z)=z^{k_d}-f_d(z)$ ($1\leq d\leq n$) are
uniformly continuous and $\|g_d\|_\infty\leq 1+a_d\leq 2$,
$\|f_{n+1}\|_\infty\leq a_{n+1}$, there exist $m_n\in\N$  and
complex numbers $c^{n,d}_j$, $b^n_j$ ($1\leq d\leq n$, $1\leq j\leq
m_n$) such that $m_n>m_{n-1}$, $m_{n-1}$ is a divisor of $m_n$  and
(P3) is satisfied. Since $\mu^n$ is absolutely continuous, there
exists $j_n\in\N$ such that $j_n>j_{n-1}$ and (P6) is satisfied.
Obviously conditions (P1) and (P2) are also satisfied.

Thus, the induction step is described and the construction of $k_n$,
$j_n$, $m_n$, $c^{n,d}_j,b^n_j$  and $\mu^n$ is complete.

First, we shall prove weak convergence of $\mu^n$ to a measure
$\mu\in\P$. Let $f\in {\cal C}(\T)$ and $\epsilon>0$. Since
$m_n\to\infty$ as $n\to\infty$ and $f$ is uniformly continuous,
there exist $a\in\N$ and complex numbers $z_1,\dots,z_{m_a}$ such
that $|f(z)-z_j|\leq\epsilon$ if $z\in I_j^a$ and $1\leq j\leq m_a$.
From (P2) and (P4) it follows that $\mu^n(I_j^a)=\mu^m(I_j^a)$ if
$m\geq a$, $n\geq a$ and $1\leq j\leq m_a$. The same argument as in
the proof of Lemma~4.1 shows that $[\mu^n-\mu^m,f]\leq 2\epsilon$
for any $m,n\geq a$. Hence $\mu^n$ is a Cauchy sequence with respect
to $\sigma$. Since, according to the Prokhorov theorem \cite{bil},
$\P$ is compact in $(\M,\sigma)$, there exists $\mu\in\P$ such that
$\mu^n\tos \mu$ as $n\to\infty$.

Next, we shall show that $\mu$ together with the sequence $k_n$
satisfy the statement of the Lemma. First, we prove that
$\mu\in\M_0$. Let $n\in\N$ and $l\in\Z$ be such that $j_n\leq
|l|<j_{n+1}$. By (P6) $|\widehat{\mu^n}(l)|\leq \epsilon_n$.
According to (P7) we have
$|\widehat{\mu^{k+1}}(l)-\widehat{\mu^{k}}(l)|\leq \epsilon_{k+1}$
for $k\geq n+1$. Finally, (P8) implies that
$|\widehat{\mu^{n+1}}(l)-\widehat{\mu^{n}}(l)|\leq 2 a_{n+2}$. Thus,
$$
|\mmu(l)|\leq |\widehat{\mu^n}(l)|+
|\widehat{\mu^{n+1}}(l)-\widehat{\mu^{n}}(l)|+ \sum_{k=n+1}^\infty
|\widehat{\mu^{k+1}}(l)-\widehat{\mu^{k}}(l)| \leq
2a_{n+2}+\sum_{j=n}^\infty \epsilon_j\to 0
$$
as $n\to\infty$. Therefore $\mmu(l)\to 0$ as $|l|\to \infty$, that
is $\mu\in\M_0$.

It remains to estimate $[\mu,g_n\overline{g_d}]$. Let $d,n\in\N$ and
$d<n$. Denote $h_n(z)=z^{k_n}$. Since $g_l=h_l-f_l$, we can write
$$
g_n\overline{g_d}=(h_n-b_j^{n-1})\overline{c_j^{n-1,d}}+(b_j^{n-1}-f_n)
\overline{c_j^{n-1,d}}+g_n(\overline{g_d}-\overline{c_j^{n-1,d}}).
$$
Taking into account that $\mu^n=\sum\limits_{j=1}^{m_{n-1}}\mu^n_j$,
where $\mu^n_j$ is the restriction of $\mu^n$ to $I^{m_{n-1}}_{j}$,
we obtain
\begin{align*}
[\mu^n,g_n\overline{g_d}]&=\sum_{j=1}^{m_{n-1}}
[\mu^n_j,g_n\overline{g_d}]= \sum_{j=1}^{m_{n-1}}
(A_j^n+B_j^n+C_j^n),\ \ \text{where}
\\
A_j^n&=[\mu^n_j,(h_n-b^{n-1}_j)\overline{c_j^{n,d}}], \
B_j^n=[\mu^n_j,(b^{n-1}_j-f_n)\overline{c_j^{n,d}}] \ \text{and}\
C_j^n=[\mu^n_j,g_n(\overline{g_d}-\overline{c_j^{n,d}})].
\end{align*}
Using that $\|g_n\|_\infty\leq 1+a_n\leq 2$ and (P3), we find that
$|c_j^{n,d}|\leq 2$ and $|b^{n-1}_j-f_n(z)|\leq \epsilon_n$ and
$|\overline{g_d}(z)-\overline{c_j^{n,d}}|\leq\epsilon_n$ for $z\in
I^{m_{n-1}}_j$.  Since $\mu^n_j$ is supported on $I^{m_{n-1}}_j$, we
have
$$
|B_j^n|\leq 2\epsilon_n\mu^n(I^{m_{n-1}}_j) \ \ \text{and}\  \
|C_j^n|\leq 2\epsilon_n\mu^n(I^{m_{n-1}}_j).
$$
On the other hand
$$
A_j^n=\overline{c_j^{n,d}}([\mu^n_j,h_n]-
b^{n-1}_j\mu^n(I^{m_{n-1}}_j))=\overline{c_j^{n,d}}
(\widehat{\mu^n_j} (k_n)-b_j^{n-1}\mu^n(I^{m_{n-1}}_j)).
$$
Using (P5) and the fact that $|c_j^{n,d}|\leq 2$, we obtain
$$
|A_j^n|\leq 2\epsilon_n/m_{n-1}.
$$
Upon putting the estimates on $A_j^n$, $B_j^n$ and $C_j^n$ together,
we have
$$
|[\mu^n,g_n\overline{g_d}]|\leq \sum_{j=1}^{m_{n-1}}
(|A_j^n|+|B_j^n|+|C_j^n|)\leq 4\epsilon_n\left(
\sum_{j=1}^{m_{n-1}} \mu^n(I^{m_{n-1}}_j)\right)+
2m_{n-1}\epsilon_n/m_{n-1}=6\epsilon_n.
$$
From (P.9) for $m\geq n>d$ it follows that
$$
|[\mu^{m+1}-\mu^m,g_n\overline{g_d}]|\leq \epsilon_{m+1}.
$$
Therefore, from (\ref{epde}) we see that
$$
|[\mu,g_n\overline{g_d}]|\leq |[\mu^n,g_n\overline{g_d}]|+
\sum_{m=n}^\infty|[\mu^{m+1}-\mu^m,g_n\overline{g_d}]|\leq
6\epsilon_n+\sum_{m=n}^\infty\epsilon_{m+1}<6\sum_{m=n}^\infty
\epsilon_{m}<\delta_n.
$$
Thus, $\mu$ satisfies all required conditions. \square

\subsection{Proof of Theorem~1.2}

We choose a set $S=\{h_n:n\in\N\}$ dense in the unit sphere of the
Banach space ${\cal C}(\T)$ and a one-to-one map
$\phi=(\phi_1,\phi_2)$ from $\N$ onto $\N^2$. We consider
$$
f_n=2^{-\phi_1(n)}(\phi_2(n))^{-1/2}h_{\phi_1(n)},
$$
which are in  ${\cal C}(\T)$, $\|f_n\|_\infty<1$ for any $n\in\N$
and $\|f_n\|_\infty\to0$ as $n\to\infty$. By Lemma~4.4 there exists
$\mu\in\P\cap\M_0$ and a strictly increasing sequence $k_n$ of
positive integers, such that
$$
|[\mu,g_n\overline{g_m}]|\leq 2^{-n}\ \ \text{whenever $n>m$},
$$
where $g_n(z)=z^{k_n}-f_n(z)$. Since $S$ is norm-dense in the unit
sphere of ${\cal C}(\T)$, we find that $\Omega=\{\lambda x:
\lambda\in\C,\ x\in S\}$ is norm-dense in ${\cal C}(\T)$, which is
in turn norm dense in $L_2(\mu)$. It follows that $\Omega$ is weakly
dense in $L_2(\mu)$,

Since $|[\mu,g_n\overline{g_m}]|=\langle g_n,g_m\rangle$, where
$\langle\cdot,\cdot\rangle$ stands for the inner product in
$L_2(\mu)$, we have
$$
\sum_{m,n\in\N\atop n>m}|\langle g_n,g_m\rangle|^2\leq
\sum_{n=2}^\infty (n-1)4^{-n}<\infty.
$$
According to Lemma~2.6, $\{g_n\}$ is a 2-sequence in $L_2(\mu)$.

We shall show that the constant function $u(z)\equiv 1$ is a weakly
supercyclic vector for the operator $Mf(z)=zf(z)$ acting on
$L_2(\mu)$. For $n\in\N$, let $A_n=\{m\in\N:\phi_1(m)=n\}$. Since
$\phi$ is one-to-one from $\N$ onto $\N^2$, it follows that $\phi_2$
is one-to-one from $A_n$ onto $\N$. Let $m\in A_n$. We have
$f_m=2^{-n}(\phi_2(m))^{-1/2}h_{n}$ and $g_m=T^{k_m}u-f_m$. Hence
$$
g_m=\beta(m)T^{\gamma(m)}-\alpha(m)h_n\ \ \text{for each}\  \ m\in
A_n,
$$
where $\beta(m)=1$, $\gamma(m)=k_m$ and $\alpha(m)=
2^{-n}(\phi_2(m))^{-1/2}$. Since $\phi_2$ is one-to-one from $A_n$
onto $\N$, we obtain
$$
\sum_{m\in A_n}|\alpha(m)|^2=
\sum_{m\in A_n}\left(2^n(\phi_2(m))^{1/2}\right)^{-2}=
2^{-2n}\sum_{j=1}^\infty j^{-1}=\infty.
$$
Upon applying Proposition~2.4, we see that $u$ is a weakly
supercyclic vector for $M$. Clearly the requirement
$\mu\in\P\cap\M_0$ also holds. The proof of Theorem~1.2 is complete.

\section{Proof of Theorem 1.6}

We start with reformulating the Salas criteria of hypercyclicity and
supercyclicity of bilateral weighted shift in a more convenient
form. This form is reminiscent of the one of Feldman \cite{feld},
which he obtained under the additional assumption of invertibility.

\prop{5.1}Let $T$ be a bilateral weighted shift acting on
$\ell_p(\Z)$ with $1\leq p<\infty$ or $c_0(\Z)$. Then $T$ is
hypercyclic if and only if for any $k\in\zp$,
\begin{equation}
\ilim\limits_{n\to\infty} \max\{\beta(k-n+1,k),
(\beta(k+1,k+n))^{-1}\}=0 \label{sal3}
\end{equation}
and $T$ is supercyclic if and only if for any $k\in\zp$,
\begin{equation}
\ilim\limits_{n\to+\infty} \beta(k-n+1,k)\beta(k+1,k+n)^{-1}=0,
\label{sal4}
\end{equation}
where $\beta(a,b)$ are the numbers defined in $(\ref{beta})$.\epr

\proof Obviously, if (\ref{sal1}) is satisfied for any $k\in\zp$
then (\ref{sal3}) holds true for any $k\in\zp$ and if (\ref{sal2})
is satisfied for any $k\in\zp$ then (\ref{sal4}) holds true for any
$k\in\zp$. It remains to prove the opposite. For any $m\in\N$ denote
$$
d_m=\bigl(\max\{1,\|w\|_\infty\}\bigr)^{-2m}\min\limits_{-m\leq
a\leq b\leq m}\beta(a,b).
$$

Suppose that (\ref{sal3}) holds for any $k\in\zp$ and (\ref{sal1})
fails for $k=m-1\in\zp$. Then there exist sequences
$\{j_n\}_{n\in\zp}$, $\{k_n\}_{n\in\zp}$ and $c>0$ such that
$|j_n|<m$, $|k_n|<m$ and $\max\{\beta(j_n-n,j_n),
(\beta(k_n,k_n+n))^{-1}\}\geq c$ for each $n\in\zp$. Since
\begin{align}
\beta(m-n+1,m)&=\frac{\beta(j_n-n,j_n)\beta(j_n+1,m)}{\beta(j_n-n,m-n)}
\geq d_m\beta(j_n-n,j_n)\ \ \text{and}\label{sal5}
\\
\beta(m+1,m+n)&\leq\frac{\beta(k_n,k_n+n)\beta(k_n+n+1,m+n)}{\beta(k_n,m)}
\leq \frac{\beta(k_n,k_n+n)}{d_m}, \label{sal6}
\end{align}
we obtain that
$$
\max\{\beta(m-n+1,m), (\beta(m+1,m+n))^{-1}\}\geq
d_m\max\{\beta(j_n-n,j_n), (\beta(k_n,k_n+n))^{-1}\}\geq cd_m
$$
for each $n\in\zp$. Thus, (\ref{sal3}) fails for $k=m$. A
contradiction.

Finally suppose that (\ref{sal4}) holds for any $k\in\zp$ and
(\ref{sal2}) fails for $k=m-1\in\zp$. Then there exist sequences
$\{j_n\}_{n\in\zp}$, $\{k_n\}_{n\in\zp}$ and $c>0$ such that
$|j_n|<m$, $|k_n|<m$ and $\frac{\beta(j_n-n,j_n)}{
\beta(k_n,k_n+n)}\geq c$ for each $n\in\zp$. Applying (\ref{sal5})
and (\ref{sal6}), we obtain
$$
\beta(m-n,m)\beta(m,m+n)^{-1}\geq d_m^2
\beta(j_n-n,j_n)\beta(k_n,k_n+n)^{-1}\geq d_m^2c
$$
for any $n\in\zp$. Thus, (\ref{sal4}) fails for $k=m$. A
contradiction. \square

The proof is based on the following two propositions on weak
closeness of sequences in $\ell_p$ with rapidly increasing norms.

\prop{5.2}Let $\H$ be a real or complex Hilbert space and
$\{x_n\}_{n\in\zp}$ be a sequence of elements of $\H$, such that
\begin{equation}
\sum_{n=0}^\infty\|x_n\|^{-a}<\infty \label{kball}
\end{equation}
for $a=2$. Then the set $S=\{x_n:n\in\zp\}$ is weakly closed in
$\H$. \epr

\prop{5.3}Let $1<p<\infty$ and $\{x_n\}_{n\in\zp}$ be a sequence of
elements of the real or complex Banach space $\ell_p(\Lambda)$, such
that $(\ref{kball})$ is satisfied with $0<a<\min\{2,q\}$, where
$\frac1p+\frac1q=1$. Then the set $S=\{x_n:n\in\zp\}$ is weakly
closed in $\ell_p(\Lambda)$. \epr

We shall now prove Theorems~1.6 with the help of these results,
postponing the proofs of the propositions to the end of the section.
For sake of completeness we formulate an analog of Propositions~5.2
and 5.3 for general Banach spaces, which we also prove in the end of
the section.

\prop{5.4}Let $\B$ be a real or complex Banach space and
$\{x_n\}_{n\in\zp}$ be a sequence of elements of $\B$, such that
$(\ref{kball})$ is satisfied with $a=1$. Then the set
$S=\{x_n:n\in\zp\}$ is weakly closed in $\B$. \epr

Note that weak closeness of a countable subset $\{x_n:n\in\zp\}$ of
a Banach space under the condition that $\|x_n\|$ grow exponentially
was proved in \cite{dt}.

\subsection{Proof of Theorem 1.6: the supercyclicity case}

In this section $\K$ stands for either the field $\R$ of real
numbers or the field $\C$ of complex numbers.

\lemma{5.5}Let $\{x_n\}_{n\in\zp}$ be a sequence in a Banach space
$\B$ over the field $\K$, $y\in\B$, $z\in\B^*$ be such that $\langle
y,z\rangle=1$ and $\Omega=\{\lambda x_n:\lambda\in\K,\ n\in\zp\}$.
If $y$ belongs to the weak closure of $\Omega$, then it belongs to
the weak closure of
$$
N = \Bigl\{\frac{x_n}{\langle x_n,z\rangle}:n\in\zp,\ \ \langle
x_n,z\rangle\neq0 \Bigr\}.
$$
\rm

\proof Let $\B_0=\{u\in \B:\langle u,z\rangle=0\}$ and consider
 $M=\B\setminus\B_0$, $\Omega_0=\Omega\cap\B_0$ and
 $\Omega_1=\Omega\setminus \B_0$. Clearly,
 $\Omega=\Omega_0\cup\Omega_1$ and $y$ is not in the weak
closure of $\Omega_0$, since $\Omega_0$ is contained in the weakly
closed set $\B_0$ and $y\notin\B_0$. Hence, $y$ is in the weak
closure of $\Omega_1$. Since the map
$$
F:M\to\B,\quad F(u)=\frac{u}{\langle u,z\rangle}
$$
is weak-to-weak continuous and $y$ is in the weak closure of
$\Omega_1$, we obtain that $F(y)=y$ is in the weak closure of
$F(\Omega_1)=N$, as required. \square

We start with a general condition for an operator to be not weakly
supercyclic.

\theorem{5.6}Let $T$ be a bounded linear operator acting on an
infinite dimensional Banach space $\B$ and $f\in\B$ be such that
$T^nf\neq 0$ for each $n\in\zp$. Assume that there exists
$y\in\B^*$, $y\neq 0$ and $a>0$ for which
\begin{equation} \sum_{n=0}^\infty\left(\frac{|\langle
T^nf,y\rangle|}{\|T^nf\|}\right)^{a}<\infty. \label{nws}
\end{equation}
Suppose also that either $a=1$ or $\B$ is a Hilbert space and $a=2$
or $\B$ is isomorphic to $\ell_p$ with $1<p<\infty$ and
$a<\min\{2,q\}$, where $\frac1p+\frac1q=1$. Then $f$ is not a weakly
supercyclic vector for $T$. \epr

\proof Since $\B$ is infinite dimensional, we can pick
$x\in\B\setminus O_{\rm pr}(T,f)$ such that $\langle x,y\rangle=1$.
Suppose that $f$ is a weakly supercyclic vector for $T$. Then $x$ is
in the weak closure of $O_{\rm pr}(T,f)$. By Lemma~5.5 $x$ is in the
weak closure of the set
$$
N=\{u_n:n\in A\}, \ \ \text{where}\ \ A=\{n\in\zp:\langle
T^nf,y\rangle\neq 0\}\ \ \text{and}\ \ u_n=\frac{T^nf}{\langle
T^nf,y\rangle}.
$$
From (\ref{nws}) it follows that $\sum\limits_{n\in
A}\|u_n\|^{-a}<\infty$. Applying Propositions~5.2, 5.3 and 5.4, we
see that $N$ is weakly closed in $\B$. Hence the only way for $x$ to
be in the weak closure of $N$ is to coincide with one of the
$u_n$'s. In this case $x\in O_{\rm pr}(T,f)$. A contradiction.
\square

Now we are ready to prove the supercyclicity part of Theorem 1.6. We
have to demonstrate that any weakly supercyclic bilateral weighted
shift acting on $\ell_p(\Z)$ with $1\leq p\leq 2$ is supercyclic.
Since according to Theorem~S, supercyclicity of a bilateral weighted
shift does not depend on $p$, by comparison principle it suffices to
consider the case $p=2$. Suppose that $T$ is a bilateral weighted
shift acting on $\ell_2(\Z)$, which is weakly supercyclic and
non-supercyclic, $w$ is its weight sequence and $\beta(a,b)$ are the
numbers defined in (\ref{beta}). By Proposition~5.1 there are $c>0$
and $m\in\zp$ such that
\begin{equation}
\beta(m-n+1,m)\geq c\,\beta(m+1,m+n)\text{\  \ for each
$n\in\zp$}.\label{salsal}
\end{equation}
Since the set of weakly supercyclic vectors of a weakly supercyclic
operator is weakly dense, there exists a weakly supercyclic vector
$x$ for $T$ in $\ell_2(\Z)$ such that $\langle x,e_m\rangle \neq0$.
Using (\ref{salsal}), we have
$$
\frac{|\langle T^nx,e_m\rangle|}{\|T^nx\|_2}\leq \frac{|\langle
x,T^{*n}e_m\rangle|}{|\langle x,e_m\rangle|\|T^ne_m\|_2}=
\frac{|\langle x,e_{n+m}\rangle|\beta(m+1,m+n)}{|\langle
x,e_m\rangle|\beta(m-n+1,m)}\leq\frac{|\langle
x,e_{n+m}\rangle|}{c|\langle x,e_m\rangle|}.
$$
Since $x\in\ell_2(\Z)$, we see that $
\sum\limits_{n=0}^\infty\left(\frac{|\langle
T^nx,e_m\rangle|}{\|T^nx\|_2}\right)^2<\infty$. By Theorem~5.6 $x$
can not be a weakly supercyclic vector for $T$. A contradiction. The
proof is complete.

\subsection{Proof of Theorem 1.6: the hypercyclicity case}

We start with the following lemma dealing with positive infinite
matrices.

\lemma{5.7}Let $\Lambda$ be an infinite countable set and
$\{a_{\alpha,\beta}\}_{\alpha,\beta\in\Lambda}$ be an infinite
matrix with non-negative entries such that
$\max\{a_{\alpha,\beta},a_{\beta,\alpha}\}\geq1$ for each
$\alpha,\beta\in\Lambda$. Then
$$
\sum_{\alpha\in\Lambda}S_\alpha^{-r}<\infty\ \ \text{for each $r>1$,
where}\ \ S_\alpha=\sum\limits_{\beta\in\Lambda}
a_{\alpha,\beta}\in[0,\infty].
$$\epr

\proof Let $M>0$ and $\alpha_1,\dots,\alpha_n$ be pairwise different
elements of $\Lambda$ such that $S_{\alpha_j}\leq M$ for $1\leq
j\leq n$. Then
$$
Mn\geq\sum_{j=1}^n S_{\alpha_j}=\sum_{1\leq j\leq n\atop
\beta\in\Lambda}a_{\alpha_j,\beta}\geq \sum_{1\leq j,k\leq
n}a_{\alpha_j,\alpha_k}=\frac12\sum_{1\leq j,k\leq
n}(a_{\alpha_j,\alpha_k}+a_{\alpha_k,\alpha_j}).
$$
Since $a_{\alpha,\beta}+a_{\beta,\alpha}\geq
\max\{a_{\alpha,\beta},a_{\beta,\alpha}\}\geq1$ for each
$\alpha,\beta\in\Lambda$, we obtain $Mn\geq n^2/2$. Hence $n\leq
2M$. Therefore for any $M>0$ there exists at most $[2M]$ elements
$\alpha$ of $\Lambda$ for which $S_\alpha\leq M$, where $[t]$ stands
for the integer part of $t\in\R$. It follows that there exists a
bijection $\phi:\N\to\Lambda$ such that the sequence $S_{\phi(n)}$
is monotonically non-decreasing. Using the above estimate with
$M=S_{\phi(n)}$, we obtain that $S_{\phi(n)}\geq n/2$ for each
$n\in\N$. Hence
$$
\sum_{\alpha\in\Lambda}S_\alpha^{-r}=\sum_{n=1}^\infty
S_{\phi(n)}^{-r}\leq \sum_{n=1}^\infty (n/2)^{-r}<\infty\ \ \text{if
$r>1$.\quad\square}
$$

Now we are ready to prove the hypercyclicity part of Theorem~1.6. We
have to demonstrate that any weakly hypercyclic bilateral weighted
shift actin on $\ell_p(\Z)$ with $1\leq p<2$ is hypercyclic. Since
according to Theorem~S, hypercyclicity of a bilateral weighted shift
does not depend on $p$, by comparison principle it suffices to
consider the case $1<p<2$. Suppose that $T$ is a non-hypercyclic
weakly hypercyclic bilateral weighted shift acting on $\ell_p(\Z)$
with $1<p<2$, $w$ is its weight sequence and $\beta(a,b)$ are the
numbers defined in (\ref{beta}). By Proposition~5.1 there are
$c\in(0,1]$ and $m\in\zp$ such that
\begin{equation}
\max\{\beta(m-n+1,m),\beta(m+1,m+n)^{-1}\}\geq c\text{\ \ for each
$n\in\zp$}.\label{salsa}
\end{equation}
Let $x$  be a weakly hypercyclic vector for $T$ and
$$
A=\{k\in\N:|\langle T^kx,e_m\rangle|>1\}.
$$
The set $\{T^kx:k\in A\}$ can not be weakly closed. Indeed,
otherwise $O(T,x)$ can not be weakly dense in the non-empty weakly
open set $\{u\in\ell_p(\Z):|\langle u,e_m\rangle|>1\}$. By
Proposition~5.3,
\begin{equation}
\sum_{k\in A}\|T^kx\|_p^{-a}=\infty\ \ \text{for each $a<2$.}
\label{53}
\end{equation}
By definition of the set $A$, we have $|\langle
x,e_{k+m}\rangle|\beta(m+1,k+m)>1$ for any $k\in A$. Hence
\begin{equation}
|\langle x,e_{k+m}\rangle|>\beta(m+1,m+k)^{-1} \ \ \text{for each
$k\in A$.} \label{531}
\end{equation}
Let now $j\in A$. Obviously
$$
\|T^jx\|_p^p=\sum_{n\in\Z}\beta(n-j+1,n)^p|\langle
x,e_n\rangle|^p\geq \sum_{k\in A}\beta(m+k-j+1,m+k)^p|\langle
x,e_{m+k}\rangle|^p.
$$
Using (\ref{531}), we obtain
\begin{align}
&\|T^jx\|_p^p\geq \sum_{k\in
A}\frac{\beta(m+k-j+1,m+k)^p}{\beta(m+1,m+k)^p}=c^p\sum_{k\in
A}a_{j,k},\label{ajk}
\\
&\quad\text{where}\ \
a_{j,k}=\left\{\begin{array}{ll}c^{-p}&\text{if $k=j$;}
\\ c^{-p}\beta(m+k-j+1,m)^p&\text{if $k<j$;}
\\ c^{-p}\beta(m+1,m+k-j)^{-p}&\text{if $k>j$}.\end{array}\right.
\notag
\end{align}
From (\ref{salsa}) it follows that
$\max\{a_{j,k},a_{k,j}\}\geq 1$ for each $j,k\in A$. Lemma~5.7
together with (\ref{ajk}) implies that
$$
\sum_{j\in A}\|T^jx\|_p^{-rp}<\infty \ \ \text{for each $r>1$}.
$$
Since $p<2$, we can choose $r>1$ such that $rp<2$. Hence the last
display contradicts (\ref{53}). The proof is complete.

\subsection{Proof of Propositions 5.3 and 5.5}

We need the following interesting theorems by Ball \cite{ball,
ball1}.

\theorem{B1}Let $\H$ be a complex Hilbert space, $\{x_n\}_{n\in\zp}$
be a sequence of elements of $\H$ such that $\|x_n\|=1$ for any
$n\in\zp$ and $\{s_n\}_{n\in\zp}$ be a sequence of positive numbers
such that $\sum\limits_{n=0}^\infty s_n^2=1$. Then there exists
$y\in\H$ such that $|\langle x_n,y\rangle|\geq s_n$ for each
$n\in\zp$. \epr

\theorem{B2}Let $\B$ be a real Banach space, $\{x_n\}_{n\in\zp}$ be
a sequence of elements of $\B$ such that $\|x_n\|=1$ for any
$n\in\zp$ and $\{s_n\}_{n\in\zp}$ be a sequence of positive numbers
such that $\sum\limits_{n=0}^\infty s_n<1$. Then there exists
$y\in\B^*$ such that $|\langle x_n,y\rangle|\geq s_n$ for each
$n\in\zp$. \epr

The real and complex versions of Propositions~5.2 and~5.4 are
equivalent to each other. Indeed the real case reduces to the
complex one by replacing the space with its complexification and the
complex case reduces to the real one just by considering the complex
space as real. Thus, it suffices to prove Proposition~5.2 in the
complex case and Proposition~5.3 in the real case. Let either $\B$
be a real Banach space or $\B=\H$ be a complex Hilbert space. Let
$y\in\B\setminus S$ and $y_n=x_n-y$, $s_n=\|y_n\|^{-1}$ for
$n\in\zp$. In the Banach space case from (\ref{kball}) with $a=1$ it
follows that $\sum\limits_{n=0}^\infty s_n=C/2<\infty$. In the
Hilbert space case from (\ref{kball}) with $a=2$ it follows that
$\sum\limits_{n=0}^\infty s_n^2=C^2<\infty$. Applying Theorem~B2 in
the Banach space case and Theorem~B1 in the Hilbert space case, we
obtain that there exists $u\in\B^*$ with $\|u\|=1$ such that
$|\langle y_n/\|y_n\|,u\rangle|\geq s_n/C$ for each $n\in\zp$.
Hence, $|\langle y_n,u\rangle|\geq C^{-1}$ for each $n\in\zp$. It
means that zero is not in the weak closure of $\{y_n:n\in\zp\}$, or
equivalently, $y$ is not in the weak closure of $S$. Since $y$ is an
arbitrary point in $\B\setminus S$, we see that $S$ is weakly
closed.

\subsection{Proof of Proposition 5.3}

The ideal way to prove Proposition 5.3 would be to use an analog of
Ball's theorem for $\ell_p$-spaces. Unfortunately, it remains
undiscovered. We use probabilistic approach to prove
Proposition~5.3.

Recall few definitions. Let $\B$ be a real Banach space and $\cal F$
be the set of linearly independent finite subsets
$Y=\{y_1,\dots,y_n\}$ of $\B^*$. Let ${\cal R}_Y$ denote the family
of  sets of the form
$$
\{x\in\B: (\langle x,y_1\rangle,\dots,\langle x,y_n\rangle) \in
B\},\ \text{where $B$ is a Borel subset of $\R^n$}.
$$
Obviously, ${\cal R}_Y$ is a sub-sigma-algebra of the Borel
sigma-algebra of $\B$. A cylindric set is any element of
$$
{\cal R}(\B) = \bigcup_{Y \in {\cal F}} {\cal R}_Y.
$$
Note that ${\cal R}(\B)$ is an algebra of subsets of $\B$, but not a
sigma-algebra if $\B$ is infinite dimensional. A cylindrical measure
on $\B$ is a finite finitely-additive, non-negative measure $\mu$ on
the algebra ${\cal R}(\B)$ such that for each $Y$ in $\cal F$, the
restriction $\mu\bigr|_{{\cal R}_Y}$ is sigma-additive. The Fourier
transform of $\mu$ is the function $\widehat\mu:\B^*\to\C $ defined
by
$$
\widehat\mu(y)= \int\limits_{\B} e^{-i\langle x,y\rangle}\,d\mu(x).
$$
This integral is with respect to a sigma-additive measure, since the
function $x\mapsto e^{-i\langle x,y\rangle}$ is bounded and ${\cal
R}_{\{y\}}$-measurable and the restriction $\mu\bigr|_{{\cal
R}_{\{y\}}}$ is sigma-additive. A cylindrical measure $\mu$ is
called gaussian if for any $Y\in\cal F$, the Borel measure
$$
\ssub{\mu}{Y}(B)=\mu(\{x\in\B: (\langle x,y_1\rangle,\dots,\langle
x,y_n\rangle) \in B\})
$$
on $\R^n$ is gaussian.

Let ${\cal S}(\B)$ be the set of bounded linear operators
$T:\B^*\to\B$ satisfying the conditions
\begin{align}
\langle Tx,y\rangle&=\langle x,Ty\rangle\ \ \text{for each
$x,y\in\B^*$},\label{sa}
\\
\langle Tx,x\rangle&\geq 0\ \ \text{for each $x\in\B^*$.}
\label{pos}
\end{align}
It is well-known, see for instance \cite{bax}, Corollary~1.2,
p.~901, that for any $T\in{\cal S}(\B)$ there exists a unique
Gaussian cylindrical measure $\mut$ on $\B$ such that
$\widehat{\mut}(x)=e^{-\frac12\langle Tx,x\rangle}$ for any
$x\in\B^*$. In this case the operator $T$ is called the {\it
covariance operator} of $\mu$. We need the following
characterization of $\sigma$-additivity of Gaussian measures on
$\ell_p$. The following theorem can be found in \cite{vah}.

\theorem{V}Let $1\leq p<\infty$ and $\mu$ be a gaussian cylindrical
measure on the real Banach space $\ell_p(\Lambda)$. Then $\mu$ is
$\sigma$-additive if and only if
\begin{align*}
\sum_{\alpha\in\Lambda} |m_\alpha|^p<\infty\ \ \text{and}\ \
\sum_{\alpha\in\Lambda} s_\alpha^{p/2}<\infty,\ \ \text{where}
\\
m_\alpha=\int\limits_{\ell_p(\Lambda)}\langle
x,e_\alpha\rangle\,d\mu(x)\ \ \text{and}\ \
s_\alpha=\int\limits_{\ell_p(\Lambda)}\langle
x,e_\alpha\rangle^2\,d\mu(x).
\end{align*}\rm

Note that finiteness of the integrals defining $s_\alpha$ imply
convergence of integrals defining $m_\alpha$. One can easily verify
that for $\mu=\mut$ with $T\in{\cal S}(\ell_p(\Lambda))$,
$m_\alpha=0$ and $s_\alpha=\langle Te_\alpha,e_\alpha\rangle$. Thus,
Theorem~V implies the following corollary.

\cor{5.8}Let $1\leq p<\infty$ and $T\in{\cal S}(\ell_p(\Lambda))$.
Then $\mut$ is $\sigma$-additive if and only if
\begin{equation*}
\sum_{\alpha\in\Lambda}\langle
Te_\alpha,e_\alpha\rangle^{p/2}<\infty.
\end{equation*}\rm

We need the following two lemmas, in which $\Lambda$ is a countable
infinite set and the spaces $\ell_p(\Lambda)$ are assumed to be {\bf
real}.

\lemma{5.9}Let $1<p,q<\infty$ be such that $\frac1p+\frac1q=1$,
$k\in\N$, $A\in{\cal S}(\ell_q(\Lambda))$ be such that
$\sum\limits_{\alpha\in\Lambda}\langle
Ae_\alpha,e_\alpha\rangle^{q/2}<\infty$ and $\{u_n\}_{n\in\zp}$ be a
sequence of vectors from $\ell_p(\Lambda)$ such that $\langle
Au_n,u_n\rangle\geq 1$ for each $n\in\zp$. Then for any sequence
$a=\{a_n\}_{n\in\zp}$ of non-negative numbers such that
$a\in\ell_k$, there exist $g_1,\dots,g_k\in\ell_q(\Lambda)$ for
which
$$
\max_{1\leq j\leq k}|\langle u_n,g_j\rangle|\geq a_n \ \ \text{for
any}\ \ n\in\zp.
$$\epr

\proof Without loss of generality we can assume that $\langle
Au_n,u_n\rangle=1$ for each $n\in\zp$. Indeed, if it is not the
case, we can replace $u_n$ by $\langle Au_n,u_n\rangle^{-1/2}u_n$.

Let $K=\{1,\dots,k\}$. For $j\in K$ and $r\in(1,\infty)$ consider
the natural projections
$P_{r,j}:\ell_r(K\times\Lambda)\to\ell_r(\Lambda)$ and natural
embeddings $J_{r_j}:\ell_r(\Lambda)\to \ell_r(K\times\Lambda)$
defined on the canonical basis as $P_{r,j}e_{l,\alpha}=e_\alpha$ and
$J_{r,j}e_\alpha=e_{j,\alpha}$. Consider the bounded linear operator
$T:\ell_p(K\times\Lambda)\to \ell_q(K\times\Lambda)$ defined by the
formula
$$
Tx=\sum_{j=1}^k J_{q,j}AP_{p,j} x.
$$
In other words $T$ is the direct sum of $k$ copies of $A$. Clearly
$T\in{\cal S}(\ell_q(K\times\Lambda))$.

Since $\langle Te_{j,\alpha},e_{j,\alpha}\rangle=\langle
Ae_\alpha,e_\alpha\rangle$ for each $(j,\alpha)\in K\times\Lambda$,
we observe that
$$
\sum\limits_{(j,\alpha)\in K\times\Lambda}\langle
Te_{j,\alpha},e_{j,\alpha}\rangle^{q/2}=k\sum_{\alpha\in\Lambda}
\langle Ae_\alpha,e_\alpha\rangle^{q/2}<\infty.
$$
By Corollary~5.8 the gaussian cylindrical measure $\mu=\mut$ on
$\ell_q(K\times\Lambda)$ is $\sigma$-additive and therefore extends
to a Borel probability measure: the measure of the entire space is 1
since the Fourier transform takes value one at zero.

Let also $u_{j,n}=J_{p,j}u_n\in\ell_p(K\times\Lambda)$, for $j\in
K$, $n\in\zp$. One can easily verify that
\begin{equation}
\text{$\langle Tu_{j,n},u_{l,n}\rangle=\delta_{j,l}$ for any $l,j\in
K$ and $n\in\zp$,} \label{norm}
\end{equation}
where $\delta_{j,l}$ is the Kronecker delta. We take $c>0$ and
consider
$$
B_{n,c}= \Bigl\{y\in\ell_q(K\times\Lambda): \sum_{j=1}^k |\langle y,
u_{j,n}\rangle |^2\leq c^2 a^2_n \Bigr\}.
$$
We shall estimate $\mu (B_{n,c})$. Consider the Borel probability
measure $\nu$ on $\R^k$ defined as
$$
\nu(B)=\mu\{y\in\ell_q(K\times\Lambda): (\langle y,
u_{1,n}\rangle,\dots,\langle y,u_{k,n}\rangle)\in B\}.
$$
From (\ref{norm}), the equality $\widehat\mu(z)=e^{-\frac12\langle
Tz,z\rangle}$ and the definition of $\nu$, it follows that the
Fourier transform of $\nu$ is $\widehat\nu(t)=e^{-|t|^2/2}$. Hence,
$\nu$ has the density $ \rho_\nu(s)= (2 \pi)^{-k/2} e^{-|s|^2/2}$.
Denote $D_b^k=\{x\in\R^k:|x|\leq b\}$. Then
$$
\mu(B_{n,c})=\nu(D^k_{ca_n})=(2\pi)^{-k/2}
\int\limits_{D^k_{ca_n}}e^{-|s|^2/2}\,ds <(2\pi)^{-
k/2}\lambda_k(D^k_{ca_n})=v_kc^{k}a_n^{k},
$$
where $\lambda_k$ is the  Lebesgue measure on $\R^k$ and
$v_k=(2\pi)^{-k/2}\lambda_k(D^k_1)$. Hence,
$$
\mu\left(\bigcup_{n=0}^\infty B_{n,c}\right)\leq\sum_{n=0}^\infty
\mu(B_{n,c})<v_kc^{k}\sum_{n=0}^\infty a^k_n.
$$
Since $a\in\ell_k$, by taking  $c$  small enough we can ensure that
$$
\mu(\Lambda_c)<1=\mu(\ell_q(K\times\Lambda)),\ \ \,\,  \text{where}\
\ \Lambda_c=\bigcup_{n=0}^\infty B_{n,c}.
$$

Therefore,  there  must be $y\in\ell_q(K\times\Lambda)\setminus
\Lambda_c$. Clearly,
$$
\sum_{j=1}^k|\langle P_{q,j}y, u_n\rangle|^2=\sum_{j=1}^k|\langle y,
u_{j,n}\rangle|^2>c^2 \,a^2_n,\quad \text{for}\ \ n\in\zp.
$$
Hence  $\max\limits_{1\leq j\leq k}|\langle g_j, u_n\rangle|>a_n$
for each $n\in\zp$, where $g_j=(\sqrt k/c)P_{q,j}y\in\ell_q(\Z)$.
\square

\lemma{5.10}Let $1<p,q<\infty$ be such that $\frac1p+\frac1q=1$, and
$\{x_n\}_{n\in\zp}$ be a sequence in $\ell_p(\Lambda)$, satisfying
{\rm (\ref{kball})} with $0<a<\min\{2,q\}$. Then there exist
$k\in\N$ and $g_1,\dots,g_k\in\ell_q(\Lambda)$ such that
$$
\max_{1\leq j\leq k}|\langle x_n,g_j\rangle|\geq 1 \ \ \text{for
any}\ \ n\in\zp.
$$\rm

\proof Denote $d=\max\{a,2a/q\}$. Since $a<\min\{2,q\}$, we see that
$d<2$ and we can choose $k\in\N$ such that $k(1-d/2)\geq a$. Let
$s_n=\|x_n\|_p^{-d}$, $a_n=\|x_n\|_p^{d/2-1}$ and
$u_n=\|x_n\|_p^{d/2-1}x_n$. From (\ref{kball}) it follows that
\begin{equation}
\sum_{n=0}^\infty s_n^r<\infty,\ \ \text{where}\ \ r=\min\{1,q/2\}
\label{sn}
\end{equation}
and that $\{a_n\}_{n\in\zp}\in\ell_k$. By Hahn--Banach theorem, for
any $n\in\zp$, we can choose $y_n\in\ell_q(\Lambda)$ such that
$\|y_n\|_q=1$ and $\langle x_n,y_n\rangle =\|x_n\|_p$. Consider the
operator
$$
A:\ell_p(\Lambda)\to\ell_q(\Lambda),\qquad Ax=\sum_{n=0}^\infty
s_n\langle x,y_n\rangle y_n.
$$
According to (\ref{sn}) the sequence $\{s_n\}$ is summable and
therefore the operator $A$ is bounded. One can easily verify that
the conditions (\ref{sa}) and (\ref{pos}) for $A$ are satisfied.
Hence $A\in{\cal S}(\ell_q(\Lambda))$. Clearly
$$
\langle Au_n,u_n\rangle=\|x_n\|_p^{d-2}\sum_{m=0}^\infty s_m\langle
x_n,y_m\rangle^2\geq s_n\|x_n\|_p^{d-2}\langle
x_n,y_n\rangle^2=\|x_n\|_p^{-d}\|x_n\|_p^{d-2}\|x_n\|_p^2=1.
$$

We shall check now that
\begin{equation}
\sum_{\alpha\in\Lambda} \langle
Ae_\alpha,e_\alpha\rangle^{q/2}<\infty. \label{q2}
\end{equation}

For any  $n\in\zp$ consider the sequence $z_n$ of non-negative
numbers with the index set $\Lambda$ defined by the formula $\langle
z_n,e_\alpha\rangle=\langle y_n,e_\alpha\rangle^2$. Let also $z$ be
the sequence defined as $\langle z,e_\alpha\rangle=\langle
Ae_\alpha,e_\alpha\rangle$. Since for any $\alpha\in\Lambda$,
$$
\langle z,e_\alpha\rangle=\langle
Ae_\alpha,e_\alpha\rangle=\sum_{n=0}^\infty s_n\langle
y_n,e_\alpha\rangle^2=\sum_{n=0}^\infty s_n\langle
z_n,e_\alpha\rangle,
$$
we see that $z=\sum\limits_{n=0}^\infty s_nz_n$ in the
coordinatewise convergence sense.

{\bf Case} $p\leq 2$. In this case $q\geq 2$. Clearly
$z_n\in\ell_{q/2}(\Lambda)$ and $\|z_n\|_{q/2}=\|y_n\|_q=1$ for each
$n\in\zp$. By (\ref{sn}) the sequence $s_n$ of positive numbers is
summable and therefore the series $\sum\limits_{n=0}^\infty s_n z_n$
is absolutely convergent in the Banach space $\ell_{q/2}(\Lambda)$.
Hence $z\in\ell_{q/2}(\Lambda)$ and (\ref{q2}) follows.

{\bf Case} $p>2$. In this case $q<2$.

Recall that for $0<\rho<1$, the space $\ell_\rho(\Lambda)$ of
sequences $x=\{x_\alpha\}_{\alpha\in\Lambda}$ in
$\ell_\infty(\Lambda) $ for which
$$
\pi_\rho(x)=\sum_{\alpha\in\Lambda }|x_{\alpha}|^\rho<\infty
$$
is no longer a Banach space. The function $\pi_\rho$ is a
pseudonorm, which turns $\ell_\rho(\Lambda)$ into a complete
metrizable topological vector space, which is not locally convex.
The pseudonorm $\pi_\rho$ satisfies the triangle inequality
$\pi_\rho(x+y)\leq \pi_\rho(x)+\pi_\rho(y)$ and the homogeneity
condition $\pi_\rho(cx)=c^\rho\pi_\rho(x)$ for $c\in\R$ and
$x,y\in\ell_\rho(\Z)$.

Clearly $z_n\in\ell_{q/2}(\Lambda)$ and
$\pi_{q/2}(t_n)=\|y_n\|^q_q=1$ for each $n\in\zp$. By (\ref{sn}), we
have $\sum\limits_{n=0}^\infty s_n^{q/2}<\infty$. From the triangle
inequality and homogeneity of $\pi_{q/2}$ it follows that the series
$\sum\limits_{n=0}^\infty s_n z_n$ is convergent in the space
$\ell_{q/2}(\Lambda)$ and therefore $z\in \ell_{q/2}(\lambda)$.
Hence (\ref{q2}) is satisfied.

Thus, in any case all conditions of Lemma~5.9 are fulfilled. Hence
there exist $g_1,\dots,g_k\in\ell_q(\Z)$ such that
$\max\limits_{1\leq j\leq k}|\langle u_n,g_j\rangle|\geq a_n$ for
any $n\in\zp$. Therefore
$$
\max_{1\leq j\leq k}|\langle
x_n,g_j\rangle|=\|x_n\|_p^{1-d/2}\max_{1\leq j\leq k}|\langle
u_n,g_j\rangle|\geq \|x_n\|_p^{1-d/2}a_n=1 \ \ \text{for any}\ \
n\in\zp.\quad \square
$$

\lemma{5.11}Let $1<p,q<\infty$ be such that $\frac1p+\frac1q=1$, and
$\{x_n\}_{n\in\zp}$ be a sequence in the real or complex Banach
space $\ell_p(\Lambda)$, satisfying $(\ref{kball})$ with
$0<a<\min\{q,2\}$. Then zero is not in the weak closure of
$\{x_n:n\in\zp\}$. \epr

\proof The real case follows immediately from Lemma~5.10. In the
complex case it suffices to notice that the complex Banach space
$\ell_p(\Lambda)$, considered as a real one, is isomorphic to the
real Banach space $\ell_p(\Lambda)$. \square

We are ready to prove Proposition~5.3. Let
$y\in\ell_p(\Lambda)\setminus S$. Applying Lemma~5.11 to the
sequence, $\{x_n-y\}_{n\in\zp}$, we see that zero is not in the weak
closure of $\{x_n-y:n\in\zp\}$. Hence $y$ is not in the weak closure
of $S$. Thus $S$ is weakly closed. The proof is complete.

\section{Concluding remarks and open problems}

We start with a few general remarks. Since the Banach space $\ell_1$
enjoys the Schur property: weak and norm convergence of sequences
are equivalent \cite{dies}, weak sequential supercyclicity and weak
sequential hypercyclicity of bounded linear operators on $\ell_1$
are equivalent to supercyclicity and hypercyclicity respectively.
For operators acting on general Banach spaces it is not true, as
follows from the example of Bayart and Matheron \cite{bm}. Next
proposition shows that it is true for operators on general Banach
spaces under the additional condition that there exists a compact
operator with dense range, commuting with the given one.

\prop{6.1} Let $T$ be a bounded linear operator acting on a Banach
space $\B$. Assume that there is a compact operator $S$, acting on
$\B$, such that $S$ has dense range and $TS=ST$. Then $T$ is weakly
sequentially supercyclic if and only if $T$ is supercyclic and $T$
is weakly sequentially hypercyclic if and only if $T$ is
hypercyclic. \epr

In order to prove Proposition~6.1 we need the following topological
lemma.

\lemma{6.2} Let $X$ and $Y$  be topological spaces and $S:X\to Y$ be
a sequentially continuous map with sequentially dense range. Let
also $A\subset X$ be a sequentially dense subset of $X$. Then $S(A)$
is sequentially dense in $Y$. \epr

\proof Let $M=[S(A)]_{\rm seq}$ be the sequential closure of $S(A)$
in $Y$. Since $S$ is sequentially continuous and $M$ is sequentially
closed, we see that $S^{-1}(M)$ is sequentially closed in $X$. Since
$A\subseteq S^{-1}(M)$ and $A$ is sequentially dense in $X$, we have
$X=S^{-1}(M)$. Hence $S(X)\subseteq M$. Since $S(X)$ is sequentially
dense in $Y$, and $M$ is sequentially closed in $Y$, we obtain
$M=Y$. \square

{\bf Proof of Propositions~6.1} \ Let $x\in \B$ be a weakly
sequentially supercyclic vector for $T$. Since $S$ is compact, it is
sequentially continuous as a map from $(\B,\sigma)$ to $(\B,\tau)$,
where $\sigma$ and $\tau$ stand for the weak topology and norm
topologies respectively \cite{rud}. Since $\tau$ is metrizable, we
have that the range of $S$ is sequentially dense in $(\B,\tau)$.
Since $O_{\rm pr}(T,x)$ is sequentially dense in $(\B,\sigma)$,
Lemma~6.2 implies that $S(O_{\rm pr}(T,x))$ is sequentially dense in
$(\B,\tau)$ and therefore norm-dense in $\B$. Taking into account
that $T$ and $S$ commute we obtain that $S(O_{\rm pr}(T,x))=O_{\rm
pr}(T,Sx)$ and therefore the projective orbit $O_{\rm pr}(T,Sx)$ is
norm dense in $\B$. Thus, $Sx$ is a supercyclic vector for $T$. The
proof of the hypercyclicity case is exactly the same. One has just
to consider the orbits instead of the projective orbits. \square

Proposition~6.1 leads to some interesting questions.

 \quest{6.3}Is it possible in Proposition~$6.1$ to replace weak
sequential supercyclicity or hypercyclicity by weak supercyclicity
or hypercyclicity? In particular, does there exist a non-supercyclic
weakly supercyclic compact operator? \epr

Bes, Chan and Sanders \cite{bcs} asked whether there exists a weakly
1-sequentially hypercyclic operator which is not norm hypercyclic.
The question remains open as well as the following ones.

\quest{6.4}Does there exist a non-hypercyclic weakly sequentially
hypercyclic operator? \rm

\quest{6.5}Does there exist a weakly sequentially hypercyclic
operator which is not weakly $1$-sequentially hypercyclic?\rm

\quest{6.6}Does there exist a weakly sequentially supercyclic
operator which is not weakly $1$-sequentially supercyclic?\epr

Finally observe that according to Proposition~1.1, Theorem~1.2
provides an example of a weakly supercyclic antisupercyclic operator
on a Hilbert space, which answers a question raised in \cite{shk}.

\subsection{Measures}

The construction of a measure in the proof of Theorem~1.2 does not
provide any control of the rate of decaying of the Fourier
coefficients. In principle it is possible to make an effective
version of the construction, but one thing is obvious: the Fourier
coefficients tend to zero extremely slowly. This motivates the
following question.

\quest{6.7}Does there exist any condition on the rate of the Fourier
coefficients $\mmu(n)$ of a Borel probability measure on $\T$
$($weaker then the trivial one: $\sum|\mmu(n)|^2<\infty)$ implying
that the multiplication operator $Mf(z)=zf(z)$ acting on $L_2(\mu)$
is not weakly supercyclic? \epr

On the other hand, it would be desirable to find simpler measures,
satisfying the assertions of Theorem~1.2.

\quest{6.8}Does there exist $\mu\in\M_0\cap\P$ being an infinite
convolution of a sequence of discrete probability measures, such
that the multiplication operator $Mf(z)=zf(z)$ acting on $L_2(\mu)$
is weakly supercyclic? What about self-similar measures? \epr

As it was remarked by Bayart and Matheron \cite{bm}, if the operator
$Mf(z)=zf(z)$ acting on $L_2(\mu)$ with $\mu\in\M_+$ is weakly
supercyclic, then $\mu$ is singular. In particular, the measure in
Theorem~1.2 is singular. It follows from the fact that if
$\mu\in\M_+$ is not singular, that is $\mu$ has a non-trivial
absolutely continuous component, then there exists $n\in\N$ such
that the operator $M^n$ is not cyclic, while the powers of any
weakly supercyclic operator are weakly supercyclic and therefore
cyclic. It is not, however, the feature of absolute continuity since
$M^3$ is not cyclic if $M$ acts on $L_2(\mu)$, where $\mu$ is the
standard Cantor measure, which is purely singular.

On the other hand if $A$ is a Borel measurable subset of $\T$ such
that $z^n\neq w^n$ for any $n\in\N$ and any different $z,w\in A$ and
$\mu\in \M_+\cap\M(A)$, then $M^n$ is cyclic for any $n\in\N$. It
follows from the observation that in this case for any $n\in\N$
there exists $\mu^n\in\M_+$ such that the operator $M^n$ acting on
$L_2(\mu)$ is unitarily equivalent to $M$ acting on $L_2(\mu^n)$.
Observe that the above property of $A$ is strictly weaker than
independence of $A$. This leads us to the following question.

\quest{6.9}Let $\mu\in\P\cap\M_0$ be such that $\supp(\mu)$ is
independent. Is $M$ acting on $L_2(\mu)$ weakly supercyclic? \epr

It worth noting that the class of measures under the hypothesis of
Question~6.9 is quite large. For instance, for any Borel measurable
set $A\subset\T$ such that the set $\P\cap\M_0\cap\M(A)$ is
non-empty, there exists a measure $\mu\in\P\cap\M_0$, whose support
is an independent subset of $A$, see \cite{gg}.

\subsection{Bilateral weighted shifts}

Theorem~1.6 together with Theorem~S characterizes weakly supercyclic
bilateral weighted shifts on $\ell_p(\Z)$ with $p\leq 2$ and weakly
hypercyclic bilateral weighted shifts on $\ell_p(\Z)$ with $p<2$.
Proposition~3.4 provides a sufficient condition of weak
supercyclicity and weak hypercyclicity of bilateral weighted shifts
on general $\ell_p(\Z)$. It is not clear whether the condition of
Proposition~3.4 is also necessary. This leads to the following
problem.

\problem{6.10}Characterize $($in terms of weight sequences$)$ weakly
supercyclic bilateral weighted shifts on $\ell_p(\Z)$ for $p>2$ and
weakly hypercyclic bilateral weighted shifts on $\ell_p(\Z)$ for
$p\geq 2$. \epr

Note also that Proposition~3.4 provides more than just a weakly
supercyclic or a weakly hypercyclic vector $x$ for a bilateral
weighted shift $T$. Namely, it ensures that $\{\lambda
T^{r_n}x:n\in\zp,\lambda\in\C\}$ or $\{T^{r_n}x:n\in\zp\}$ are
weakly dense for an exponentially growing sequence of $\{r_n\}$ of
positive integers. Indeed, condition (W2) of Proposition~3.4 implies
that $\ilim\limits_{n\to\infty}(r_n)^{1/n}\geq\frac{\sqrt5+1}{2}>1$.
One way to approach Problem~6.4 could be to find out whether there
exists a weakly supercyclic or a weakly hypercyclic bilateral
weighted shift $T$ such that the sets of the shape $\{\lambda
T^{r_n}x:n\in\zp,\lambda\in\C\}$ or $\{T^{r_n}x:n\in\zp\}$ are not
weakly dense for any exponentially growing sequence of $\{r_n\}$ of
positive integers.

Using Proposition~3.4 and the technique of the proof of Theorem~1.6
it is possible for any $p\geq 2$ to find a bilateral weighted shift,
which is weakly hypercyclic on $\ell_p(\Z)$ and not weakly
hypercyclic on $\ell_r(\Z)$ for each $r<p$. Thus, the infinum of
$p$'s for which a given bilateral weighted shift is weakly
hypercyclic on $\ell_p(\Z)$ is a parameter taking all values between
$2$ and $\infty$. Thus, any characterization of hypercyclic
bilateral weighted shifts on $\ell_p(\Z)$ for $p\geq2$ must depend
on the parameter $p$.

\subsection{Tightness of Propositions 5.2--5.4}

The following theorem is known as Dvoretzky theorem on almost
spherical sections \cite{Dv}. Somewhat weaker version of this
theorem was obtained earlier by Dvoretzky and Rogers \cite{DR}.

\theorem{D} For each $n\in\N$ and each $\epsilon>0$, there exists
$m=m(n,\epsilon)\in\N$ such that for any Banach space $\B$ with
$\text{\rm dim}\,\B\geq m$ there is an $n$-dimensional linear
subspace $L$ in $\B$ and a basis $e_1,\dots,e_n$ in $L$ for which
$$
\left\|\sum_{j=1}^nc_je_j\right\|_{\B}\leq
\left(\sum_{j=1}^n|c_j|^2\right)^{1/2}\leq (1+\epsilon)
\left\|\sum_{j=1}^nc_je_j\right\|_{\B}\ \ \text{for any
$(c_1,\dots,c_n)\in\C^n$}.
$$
\epr

We use this theorem in order to prove the following proposition,
which allows us to demonstrate tightness of Propositions~5.2--5.4.

\prop{6.11}For any infinite dimensional Banach space $\B$ and any
sequence $\{c_n\}_{n\in\zp}$ of positive numbers such that
$\sum\limits_{n=0}^\infty c_n^{-2}=\infty$, there exists a sequence
$\{x_n\}_{n\in\zp}$ in $\B$ such that $\|x_n\|=c_n$ for each
$n\in\zp$ and zero is in the weak closure of $\{x_n:n\in\zp\}$. \epr

\proof Pick a strictly increasing sequence $\{n_k\}_{k\geq 0}$ of
integers such that $n_0=0$ and
\begin{equation}
\lim_{k\to\infty}\sum_{j=n_{k-1}}^{n_k-1}c_j^{-2}=\infty. \label{1}
\end{equation}
Denote $j_k=n_k-n_{k-1}$, $k\in\N$. By Theorem~D, for each $k\in\N$,
there exist a linear subspace $F_k$ of $\B$ with
$\text{dim}\,F_k=j_k$ and a basis $e_{n_{k-1}},\dots,e_{n_k-1}$ in
$F_k$ such that
\begin{equation}
\biggl\|\sum_{j=n_{k-1}}^{n_k-1}c_je_j\biggr\|_{\B}\leq
\biggl(\sum_{j=n_{k-1}}^{n_k-1}|c_j|^2\biggr)^{1/2}\leq 2
\biggl\|\sum_{j=n_{k-1}}^{n_k-1}c_je_j\biggr\|_{\B} \label{2}
\end{equation}
for any complex numbers $c_j$. In what follows, we assume that
$F_k$'s carry the norm inherited from $\B$. The inequality for the
dual norm reads as follows
\begin{equation}
\frac12\ssub{\|f\|}{F_k^*}\leq
\biggl(\sum_{j=n_{k-1}}^{n_k-1}|\langle f,e_j\rangle|^2\biggr)^{1/2}
\leq \ssub{\|f\|}{F_k^*}\ \ \text{for each}\ \ f\in F_k^*. \label{3}
\end{equation}
Denote $x_n=c_ne_n/\|e_n\|$ for $n\in\zp$. Obviously $\|x_n\|=c_n$.
It remains to prove that zero is in the weak closure of
$\{x_n:n\in\zp\}$. Suppose the contrary. Then there exist
$g_1,\dots,g_m\in\B^*$ such that
\begin{equation}
\max_{1\leq l\leq m}|\langle g_l,x_n\rangle|\geq 1\ \ \text{for each
$n\in\zp$}. \label{4}
\end{equation}
Denote $M=\max\limits_{1\leq l\leq m}\ssub{\|g_j\|}{\B^*}$ and for
each positive integer $k$ let $h_l^k\in F_k^*$ be the restriction of
$g_l$ to $F_k$. From (\ref{2}) it follows that $\|e_n\|\geq 1/2$ for
each $n$. If $1\leq l\leq m$ and $n_{k-1}\leq n\leq n_k-1$, then
$$
|\langle h^k_l,e_n\rangle|=|\langle g_l,e_n\rangle|=
\|e_n\|c_n^{-1}|\langle g_l,x_n\rangle|\geq (2c_n)^{-1} |\langle
g_l,x_n\rangle|.
$$
Using (\ref{4}) and the last display, we obtain
$$
\sum_{l=1}^m|\langle h^k_l,e_n\rangle|^2\geq (2c_n)^{-2}\ \
\text{for}\ \ n_{k-1}\leq n\leq n_k-1.
$$
Taking (\ref{3}) into account, we get
$$
\sum_{l=1}^m\ssub{\|h_l^k\|^2\!\!\!}{F_k^*}\geq
\sum_{l=1}^m\sum_{j=n_{k-1}}^{n_k-1} |\langle
h^k_l,e_j\rangle|^2\geq \frac14 \sum_{j=n_{k-1}}^{n_k-1} c_j^{-2}.
$$
Since $\ssub{\|h_l^k\|}{F_k^\star}\leq \ssub{\|g_l\|}{\B^*}\leq M$,
we see that  $4mM^2\geq \sum\limits_{j=n_{k-1}}^{n_k-1} c_j^{-2}$
for any positive integer $k$, which contradicts (\ref{1}). \square

\cor{6.12}Let $1\leq p\leq\infty$, $\B_p=\ell_p$ if $1\leq p<\infty$
and $\B_\infty=c_0$ and $\{c_n\}_{n\in\zp}$ be a sequence of
positive numbers such that $\sum\limits_{n=0}^\infty
c_n^{-r}=\infty$, where $r=\min\{2,q\}$ and $\frac1p+\frac1q=1$.
Then there exists a sequence $\{x_n\}_{n\in\zp}$ in $\B$ such that
$\|x_n\|=c_n$ for each $n\in\zp$ and zero is in the weak closure of
$\{x_n:n\in\zp\}$. \epr

\proof The case $1\leq p\leq 2$ follows from Proposition~6.11. If
$p>2$, we can take $x_n=c_ne_n$ and apply Lemma~2.2. \square

Corollary~6.12 for $p=2$ and $p=\infty$ implies that conditions on
the growth of $\|x_n\|$ in Propositions~5.2 and 5.4 are best
possible. Proposition~5.3 and Corollary~6.12 lead to the natural
conjecture that the best possible condition on the growth of
$\|x_n\|$ implying weak closeness of $\{x_n:n\in\zp\}$ in $\ell_p$
is (\ref{kball}) with $a=\min\{2,q\}$. In order to prove this
conjecture it would suffice to answer the following question
affirmatively.

\quest{6.13}Let $1\leq p<\infty$, $\{x_n\}_{n\in\zp}$ be a sequence
in the unit sphere of a the complex Banach space $\ell_p$ and
$\{s_n\}_{n\in\zp}$ be a sequence of positive numbers such that
$\sum\limits_{n=0}^\infty s_n^r=1$, where $r=\min\{2,q\}$ and
$\frac1p+\frac1q=1$. Does there exist $y\in\ell_q$ such that
$\|y\|_q=1$ and $|\langle x_n,y\rangle|\geq s_n$ for each $n\in\zp$?
\epr

Note that an affirmative answer to this question would also provide
an interesting generalization of Ball's theorem (Theorem~B1) and
possibly lead to further applications in convex analysis.

\subsection{Sequential weak topology}

In this final section we discuss the nature of weak sequential
density and thus of weak sequential supercyclicity and
hypercyclicity. Recall that a topological space $(X,\tau)$ is called
{\it sequential} if a subset of $X$ is closed if and only if it is
sequentially closed.  A subset $A$ of a topological vector space
$(X,\tau)$ is called sequentially open if $X\setminus A$ is
sequentially closed. It is straightforward to verify that the
collection $\tau_{\rm seq}$ of sequentially open subsets of a
topological space $(X,\tau)$ forms a topology. Moreover,
$\tau\subseteq\tau_{\rm seq}$ and $(X,\tau_{\rm seq})$ is sequential
and a sequence converges in $(X,\tau)$ if and only if it converges
to the same limit in $(X,\tau_{\rm seq})$.

For a Banach space $\B$, $\sigma=\sigma(\B,\B^*)$ stands for the
weak topology of $\B$ and $\sigma_{\rm seq}$ stands for the
corresponding sequential topology. From the above it follows that a
set $A\subseteq \B$ is weakly sequentially dense in $\B$ if and only
if it is dense in $\sigma_{\rm seq}$. Thus, the concepts of weak
sequential hypercyclicity and supercyclicity (unlike weak
1-sequential hypercyclicity and supercyclicity) are topological.
Namely they are just hypercyclicity and supercyclicity with respect
to the topology $\sigma_{\rm seq}$ intermediate between the weak and
the norm topologies.

Finally we make a few remarks on the nature of the topology
$\sigma_{\rm seq}$. From the Schur Theorem \cite{dies} it follows
that the topology $\sigma_{\rm seq}$ on the Banach space $\ell_1$
coincides with the norm topology. In \cite{rest} it is observed that
there are Banach spaces $\B$ for which $(\B,\sigma_{\rm seq})$ fails
to be a topological vector space: the addition $(x,y)\mapsto x+y$,
although being separately continuous, may fail to be continuous. It
is also demonstrated in \cite{rest} that if $\B^*$ is separable then
$\sigma_{\rm seq}$ coincides with the so-called bounded weak
topology, which is the strongest topology that agrees with the weak
topology on the bounded sets. According to the Banach--Dieudonn\'e
theorem, see for instance \cite{sch}, the bounded weak topology on a
reflexive Banach space coincides with the pre-compact convergence
topology, that is the topology of uniform convergence over the norm
pre-compact subsets of $\B^*$. It worth mentioning that $\B$ with
the pre-compact convergence topology is a complete locally convex
topological vector space. For a characterization of local convexity
of the bounded weak topology we refer to \cite{gil}. Thus, we have
the following

\prop{6.14}Let $\B$ be a separable reflexive Banach space. Then the
weak sequential topology $\sigma_{\rm seq}$ on $\B$ coincides with
the pre-compact convergence topology. \epr

According to this proposition weak sequential supercyclicity and
hypercyclicity of bounded linear operators on a separable reflexive
Banach space are exactly supercyclicity and hypercyclicity with
respect to the pre-compact convergence topology. Note that for
infinite dimensional Banach spaces the pre-compact convergence
topology is strictly stronger than the weak topology and strictly
weaker than the norm topology.

\bigskip

{\bf Acknowledgements.} \ Partially supported by Plan Nacional I+D+I
Grant BFM2003-00034, Junta de Andaluc\'{\i}a FQM-260 and British
Engineering and Physical Research Council Grant GR/T25552/01. The
author is grateful to Alfonso Montes--Rodr\'iguez and Fr\'ed\'eric
Bayart for many helpful comments.

\end{document}